%% file: main.tex
\documentclass[hidelinks,onefignum,onetabnum]{siamart220329}

\input{ex_shared}

\input{math_commands}

\usepackage{amsmath,amssymb,a4wide}
\usepackage{latexsym}
\usepackage{indentfirst}
\usepackage{caption}
\usepackage{placeins}
\usepackage{enumerate}
\usepackage{booktabs}
\usepackage{algorithm}
\usepackage{algorithmic}
\usepackage{multirow}
\usepackage{optidef}
\usepackage{graphicx,color}
\usepackage{diagbox}
\usepackage[normalem]{ulem}
\usepackage{hyperref}
\usepackage{cleveref}
\usepackage{subcaption} %

\usepackage{tikz}
\input{ht_diagram}

\ifpdf
\hypersetup{
  pdftitle={Non-negative_HT},
  pdfauthor={X. Tang, H. Chen, and L. Ying}
}
\fi

\newcommand{\ket}[1]{\lvert #1 \rangle}

\usepackage{xcolor}

\begin{document}

\maketitle

\begin{abstract}
    In this work, we present an efficient method to compress a high-dimensional discrete probability function, i.e., a probability tensor, into a non-negative hierarchical Tucker format. The methodology is a two-stage procedure. In the first stage, we take an existing interpolation method to compress the target tensor into a hierarchical Tucker (HT) in a manner similar to the CUR decomposition for low-rank matrix reconstruction. In the second stage, we fit the first-stage output against a non-negative hierarchical Tucker ansatz using a second-order method tailored specifically for this setting. When the tensor is of order \(d\), both stages admit an \(\mathcal{O}(d)\) computational complexity, and therefore the proposed methodology readily extends into high-dimensional settings. Numerical experiments show success in compressing various high-dimensional probability tensors.
\end{abstract}



\section{Introduction}

This work proposes a new method to compress high-dimensional discrete distribution functions into tractable low-rank formats. Specifically, writing \([n] := \{1, \ldots, n\}\), the goal is to compress a \(d\)-dimensional reference distribution function \(P \colon [n]^d \to \R_{\geq 0}\). Storing \(P\) has an \(\mathcal{O}(n^d)\) complexity: even in simple cases, the dimension \(d\) can be on the order of tens of variables. Moreover, computing properties of \(P\) often requires exponential scaling in \(d\). Approximation of \(P\) within a parametric family of models \(\{P_{\theta}\}_{\theta \in \Theta}\), low-rank representation being one example, is a common approach for parameter reduction in high-dimensional cases. 

The task of compressing \(P\) falls under two cases.
The first case is variational inference (VI) \cite{jordan1999introduction,blei2017variational}. In that case, for an arbitrary multi-index \((i_1, \ldots, i_d) \in [n]^d\), the VI setting assumes that one can access \(P(i_1, \ldots, i_d)\) up to an unknown normalization constant.
The second case is density estimation (DE) \cite{silverman2018density}.
In this case, we assume that one has a collection of samples \(\left(y_{1}^{(j)}, \ldots, y_{d}^{(j)}\right)_{j=1}^{N} \subset [n]^d\) which are distributed according to \(P\). This work proposes an end-to-end algorithm for compressing \(P\), and we assume access to \(P\) either from the VI case or from the DE case.

\begin{figure}
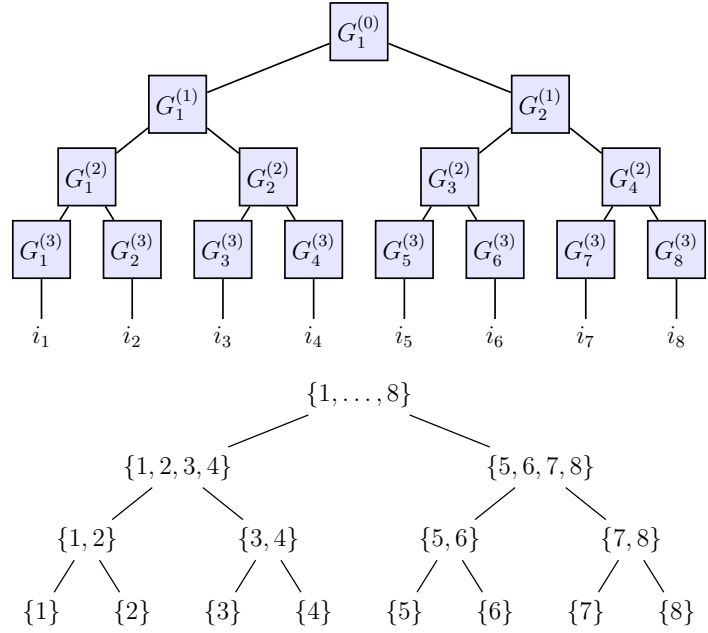

    \centering
    \scalebox{0.8}{\HTFigure}
    \caption{Hierarchical Tucker model for \(d = 8\) under a complete binary tree structure. The top figure illustrates the tensor network structure, and the bottom figure is the corresponding dimension tree, cf.\ \cite{hackbusch2009new}.}
    \label{fig: NHT}
\end{figure}

This work uses a non-negative hierarchical Tucker (NHT) format as illustrated in \Cref{fig: NHT}. For an NHT ansatz, all factors in \Cref{fig: NHT} would be entry-wise non-negative, which ensures non-negativity of the approximated tensor in full tensor space. Our end-to-end algorithm contains two stages. The first stage uses existing interpolation algorithm in the VI case and existing density estimation algorithm in the DE case to construct a hierarchical Tucker \(\tilde{P} \approx P\). The second stage uses an NHT ansatz to fit \(\tilde{P}\) by variationally optimizing the non-negative factors. The two-stage approach efficiently gives \(P\) an approximation in NHT format. Previously, there were algorithms to generate an approximation of \(\tilde{P}\) under the hierarchical Tucker format. In the DE case, works such as \cite{peng2023generative,tang2023generative,tang2024solving} have considered a hierarchical Tucker approximation of \(P\). In the VI case, one can refer to \cite{ballani2013black, ryzhakov2024black} for interpolation algorithms in constructing \(\tilde{P} \approx P\) by querying entries of \(P\). The main contribution of our work comes from using the NHT format and designing the NHT fitting procedure in the second stage.

\subsection{Background}

This section gives the background of tensor networks for readers. We give an exposition of the ansatz. Subsequently, we justify the use of the ansatz by comparison with other design alternatives such as neural network models, tensor network models without a positivity constraint, and the non-negative tensor train model.

\paragraph{Non-negativity in low-rank representation}
We illustrate the idea of a non-negative tensor network in the simplest \(d = 2\) case, where the proposal would reduce to non-negative matrix factorization (NMF) \cite{lee1999learning}. In the bivariate case, we have \(P\) being an \(n \times n\) matrix, and we compress \(P\) by finding two non-negative low-rank factors \(U, V \in \R^{n \times r}_{\geq 0}\) so that \(P = UV^{\top}\). When \(r \ll n\), the \(\mathcal{O}(nr)\) cost to store \(U, V\) is more efficient than the \(\mathcal{O}(n^2)\) cost of storing \(P\). The matrix \(UV^{\top}\) has only non-negative entries under the following simple argument: matrix multiplication involves only multiplication and addition, and \(U, V\) have only non-negative entries. 

The role of non-negativity in the low-rank factors is best illustrated by the challenges practitioners encounter when the \(U, V\) factors are unconstrained. When \(UV^{\top}\) admits at least one negative entry, one essentially has a sign problem, which leads to a wide range of undesirable consequences for downstream applications. For example, the Kullback–Leibler divergence between \(P\) and \(UV^{\top}\) is undefined in either direction. Even for the simple task of sampling from \(P\), one would find the approximation \(UV^{\top}\) hard to use, as the matrix \(UV^{\top}\) is at best only a signed measure. Therefore, non-negativity is a simple approach to retain a low-rank format while avoiding the sign problem in approximating discrete distribution functions.

\paragraph{Tensor network representations}
In what follows, we will refer to \(P \colon [n]^d \to \R_{\geq 0}\) as a \(d\)-tensor, short for a tensor of order \(d\). For the high-dimensional cases, storage and manipulation of \(P\) as a \(d\)-tensor are often achieved through a tensor network format \cite{biamonte2017tensor, white1992density, hackbusch2009new, oseledets2011tensor}. In simple terms, a tensor network stores a \(d\)-tensor by means of a collection of tensor factors of a much smaller order joined by tensor contractions. For instance, in the so-called tensor-train format \cite{oseledets2010tt,oseledets2011tensor}, known more commonly in the quantum chemistry community as matrix product state (MPS) with open boundary condition \cite{fannes1992finitely, white1992density, ostlund1995thermodynamic, vidal2003efficient}, one compresses an order-\(d\) tensor into a collection of \(d\) tensor factors of order two or three. Essentially, if one imagines the variables \{1, \ldots, d\} arranged along a link graph, then the \(d\) tensor factors consist of (i) a matrix \(G_1 \in \R^{n \times r_1}\) at the left boundary, (ii) a collection of 3-tensors  \( \{G_{k} \in \R^{r_{k-1} \times n \times r_{k}}\}_{k = 2}^{d-1}\) in the interior, and (iii) a matrix \(G_d \in \R^{r_{d-1} \times n}\) at the right boundary. When \(P\) is compressed by the tensor-train with factors \(\{G_k\}_{k =1}^{d}\), the evaluation of one entry of \(P\) reads:
\begin{equation}\label{eqn: TT}
    P(i_1, \ldots, i_d) = G_1(i_1, :)G_2(:, i_2, :)\cdots G_{d-1}(:, i_{d-1}, :)G_d(:, i_d),
\end{equation}
where \(G_1(i_1, :)\) is a row vector slice of \(G_1\), and in the interior \(G_k(:, i_k, :)\) is a matrix slice of \(G_k\), and at the right boundary \(G_d(:, i_d)\) is a column vector slice of \(G_d\). Similarly, to calculate the sum of \(P\), i.e. \(Z = \sum_{i_1, \ldots, i_d = 1}^{n}P(i_1, \ldots, i_d)\), one can adapt the formula in \Cref{eqn: TT} to obtain 
\begin{equation}\label{eqn: TT Z}
\begin{aligned}
    Z &= \sum_{i_1, \ldots, i_d = 1}^{n}P(i_1, \ldots, i_d) \\ &= \left(\sum_{i_1 = 1}^{n} G_1(i_1, :)\right)\left(\sum_{i_2 = 1}^{n}G_2(:, i_2, :)\right)\cdots \left(\sum_{i_d = 1}^{n} G_d(:, i_d)\right).
\end{aligned}
\end{equation}

It would be natural to ask why a practitioner would choose a tensor network model over alternatives. For instance, reasonable alternative approaches for density approximation include the exponential family model \cite{wainwright2008graphical, barndorff2014information, brown1986fundamentals}, the energy-based model \cite{hinton2002training, lecun2006tutorial}, and the diffusion model \cite{song2019generative,song2021maximum, lou2023discrete}, to name a few. The case for using tensor network models falls into two main arguments. The first reason is that the tensor network format, as illustrated in \Cref{eqn: TT}, is naturally defined for discrete distributions, whereas numerous neural network models are designed primarily for continuous distributions. Notably, for popular and empirically successful models such as normalizing flow \cite{tabak2010density,rezende2015variational}, flow matching \cite{lipman2022flow}, and stochastic interpolant \cite{albergo2025stochastic, albergo2023stochastic}, adapting such models to discrete distributions would require significant architectural redesign. The second reason is that, as illustrated in \Cref{eqn: TT Z}, tensor network models are generally amenable to the calculation of the normalizing constant. A calculation of \(Z\) as in \Cref{eqn: TT Z} is essential for basic tasks such as calculating the negative log-likelihood (NLL). One can see that the calculation of \(Z\) in \Cref{eqn: TT Z} only requires \(\mathcal{O}(d)\) complexity. In contrast, the calculation of \(Z\) for exponential family models and energy-based models usually requires approximations, relaxations, or other heuristics. Therefore, for discrete distribution functions, using a tensor network model allows one to compute the normalizing constant, to perform sampling, and to evaluate the log-likelihood without using heuristics. Tensor network models are thus appealing in density estimation for discrete functions. When one only needs to generate samples and the information of the normalized density is not necessary, neural network methods can be more favorable.

\paragraph{Non-negative tensor works}
The preceding exposition naturally motivates the introduction of non-negative tensor networks as a parametric family for approximating discrete distribution functions in high dimensions, i.e., large \(d\). A non-negative tensor-train (NTT) ansatz admits the formula in \Cref{eqn: TT} with the condition that every entry of \(\{G_k\}_{k = 1}^{d}\) is non-negative. As the \(d = 2\) case illustrates, when the tensor-train approximation of \(P\) has even one negative entry, generating samples or calculating the negative log-likelihood can be conceptually problematic. The NTT ansatz is a theoretically appealing approach to approximating a density. The primary benefit of a non-negative tensor network format is that the ansatz is guaranteed to be entry-wise non-negative, like energy-based models and diffusion models. The first and last authors of this work co-authored \cite{tang2025variational}, which introduces a complete end-to-end procedure to compress a discrete tensor \(P\) into an NTT format. While even the \(d = 2\) NMF case is well-known to be an NP-hard problem, the numerical result in \cite{tang2025variational} shows that compressing \(P\) into an NTT format can be done with a surprisingly high accuracy even in high-dimensional cases.

\paragraph{Non-negative hierarchical Tucker}
This work compresses \(P\) into a non-negative hierarchical Tucker (NHT) format. Essentially, an NHT format is a tensor network ansatz based on a binary tree graph, which may be a complete binary tree or an incomplete binary tree. Moreover, just like the NTT case, an NHT ansatz decomposes an order \(d\) tensor into a collection of tensor factors of order two or three. The aim of using NHT is to develop a non-negative tensor network model better suited to 2D lattice models or other general high-dimensional distributions with more intricate interaction structures. Notably, it is well-known that the tensor-train model is more suitable for compressing densities satisfying an area law \cite{eisert2008area}. Therefore, the NHT model can tackle cases where the area law assumption does not hold, such as Ginzburg-Landau models over a 2D lattice geometry \cite{ginzburg2009theory, hoffmann2012ginzburg, hohenberg2015introduction,weinan2004minimum}, and indeed the numerical experiment sections in \Cref{sec: numerics} largely consider 2D lattice problems where one cannot assume an area law to hold. 

\subsection{Main contribution}
This work can be viewed as extending the NTT compression algorithm introduced in \cite{tang2025variational} to the case of non-negative hierarchical Tucker. Similar to \cite{tang2025variational}, this work consists of two stages. In the first stage, when we are in the variational inference case, we are given a function handle of \(P\), i.e., an oracle model to query arbitrary entries of \(P\). Using existing interpolation methods based on interpolative decomposition (see details in \Cref{sec: alg}), the first stage compresses \(P\) into a hierarchical Tucker format \(\tilde{P}\). For the density estimation case, we likewise use existing methods to obtain \(\tilde{P}\) from samples of \(P\). In the second stage, one uses a variational approach to use a non-negative hierarchical Tucker ansatz \(P_{G}\) to fit against \(\tilde{P}\). Let \(G = (G_1, \ldots, G_{K})\) be the tensor factors in \(P_{G}\). Then, NHT fitting is defined over a minimization task where the loss function is 
\begin{equation}\label{eqn: NHT fitting variational form}
    \ell(G) = \lVert P_{G} - \tilde{P} \rVert_{F}^{2} + \sum_{q}\mu_q \phi_q(G_q),
\end{equation}
where \(\mu_q > 0\) is the regularization strength on \(G_q\), and \(\phi_q(G_q) = -\sum_{\alpha, \beta, \gamma}\log G_q(\alpha, \beta, \gamma)\) is a log-barrier on the entries of \(G_q\).

The second stage is the main contribution of this work. When restricted to optimizing each tensor component, one can see that \Cref{eqn: NHT fitting variational form} is a strongly convex loss function in \(G_q\). Therefore, we propose an alternating minimization approach whereby the variables \(G_q\) are updated sequentially using second-order methods. Due to the structural difference between the tensor-train model and the hierarchical Tucker model, the second stage logic is redesigned substantially. Most notably, we introduce a new warm initialization strategy. Compared to adapting the warm initialization strategy in \cite{tang2025variational} to the binary tree cases, the newly designed approach is shown to significantly address the issue of reaching local minima in high-dimensional settings, and we often see the final accuracy improving by one order of magnitude from the improved initialization. Moreover, adapting the second stage code to the case of hierarchical Tucker is a challenge in itself. Lastly, while this manuscript is illustrated over complete binary trees, our code is written to work with incomplete binary trees.

\subsection{Outline}
This work is organized as follows. \Cref{sec: alg} details the proposed NHT compression subroutine. \Cref{sec: numerics} shows the numerical performance of the proposed approach on practical variational inference and density estimation tasks. \Cref{sec: conclusion} gives concluding remarks.

\section{Main algorithm}\label{sec: alg}

This section develops the two-stage compression procedure. In the first stage, we compress the target \(P\) into a signed hierarchical Tucker reference \(\tilde{P}\). In the second stage, we fit a non-negative hierarchical Tucker ansatz \(P_{G}\) to that reference. \Cref{sec: HT} recalls the hierarchical Tucker ansatz and presents the first stage. \Cref{sec: NHT fitting} sets up the second stage as an alternating minimization procedure with a per-component Newton step. \Cref{sec: NHT fitting extra} introduces acceleration techniques. \Cref{sec: NHT warm start} presents the warm initialization algorithm at the start of the second stage. For simplicity, the main text treats the complete binary tree of depth \(L\), so that \(d = 2^{L}\). \Cref{sec: appendix TTNS} extends the method to general tree tensor networks.

\paragraph{Motivating example}
High-dimensional non-negative tensors arise from the grid discretization of an unnormalized Boltzmann distribution \(p \propto \exp(-V)\). We take the two-dimensional Ginzburg-Landau model \cite{ginzburg2009theory, hoffmann2012ginzburg, hohenberg2015introduction, weinan2004minimum} as the running example. A state is a scalar field \(x \colon [0,1]^2 \to \R\) with potential
\begin{equation}\label{eqn: GL potential continuous}
  V(x) = \frac{\lambda}{2}\int_{[0,1]^2} \lvert \nabla_{s} x(s)\rvert^2 \, ds
       + \frac{1}{4\lambda}\int_{[0,1]^2} \left(1 - x(s)^2\right)^2 \, ds,
\end{equation}
where \(\lambda\) sets the balance between the gradient term and the double-well term.
Discretizing \([0,1]^2\) on an \(m \times m\) grid of spacing \(h\) replaces the field by the vector \(x = (x_{v})_{v}\) with \(x_{v} = x(s_{v})\). The dimension is \(d = m^2\). The potential becomes
\begin{equation}\label{eqn: GL potential discrete}
  V(x) = \frac{\lambda}{2}\sum_{v \sim w}\left(\frac{x_{v} - x_{w}}{h}\right)^2
       + \frac{1}{4\lambda}\sum_{v}\left(1 - x_{v}^2\right)^2,
\end{equation}
where \(v \sim w\) denotes a pair of adjacent grid points. Restricting each \(x_{v}\) to \(n\) grid values then gives the order-\(d\) non-negative tensor \(P(i_1, \ldots, i_d) = \exp(-V(x_{i_1}, \ldots, x_{i_d}))\) that we compress. This example fits both access models. The variational inference case assumes one can read off \(P(i_1, \ldots, i_d)\) by evaluating \(V\). The density estimation case can similarly use the analytic expression of \(P\) to generate samples by performing Markov chain Monte Carlo algorithms, and we point readers to \cite{liu2001monte} for a reference. Notably, the neighbor coupling \(v \sim w\) is two-dimensional, which is less suited to the linear geometry in the tensor-train case. The binary tree structure in hierarchical Tucker is better suited for this coupling. Concretely, for the 2D lattice case, one can recursively apply bipartition to the first and second axes, and the resultant tree structure provides a good inductive bias to capture the global correlation structure of neighboring points on the 2D grid.

\paragraph{Notation}
We fix the notation for the section. For \(n \in \mathbb{N}\), we write \([n] := \{1, \ldots, n\}\). For \(S \subseteq [d]\), we use \(i_{S}\) to denote the subvector of \(i\) indexed by \(S\), and we use \(\bar{S} := [d] \setminus S\) to denote the complement. The complete binary tree has depth \(L\) and \(d = 2^{L}\) leaves. Throughout the work, we let \(q\) denote the node at level \(l\) and position \(k \in [2^{l}]\). The node \(q\) is illustrated in \Cref{fig: NHT} as \(G^{(l)}_{k}\). We use the identification \(q := I_{k}^{(l)}\). Level \(l\) partitions the variables as
\begin{equation}\label{eqn: bipartition}
  [d] = \bigcup_{k=1}^{2^{l}} I_{k}^{(l)}, \qquad I_{k}^{(l)} := \{ 2^{L-l}(k-1) + 1, \ldots, 2^{L-l}k \},
\end{equation}
and the block of \(q\) splits into the blocks of its two children, \(I_{k}^{(l)} = I_{2k-1}^{(l+1)} \cup I_{2k}^{(l+1)}\). For a fixed node \(q\), we use the identification
\[
  a := I_{2k-1}^{(l+1)}, \qquad b := I_{2k}^{(l+1)}, \qquad f := [d] \setminus I_{k}^{(l)},
\]
which respectively correspond to the block of the left child of \(q\), the block of the right child of \(q\), and the variables corresponding to the non-descendants of \(q\). We write \(a\), \(b\), \(f\) when \(q\) is understood. The reference \(\tilde{P}\) has signed tensor components \((F_{q})_{q}\) of maximal internal rank \(r\). The NHT ansatz \(P_{G}\) has non-negative tensor components \((G_{q})_{q}\) of maximal internal rank \(\rho\).

\subsection{Stage one: hierarchical Tucker interpolation}\label{sec: HT}

The first stage compresses \(P\) into a hierarchical Tucker tensor, and we cover the interpolation setting for variational inference. The output \(\tilde{P}\) of this stage has signed components in general and can not guarantee non-negativity. For the density estimation case, we refer the readers to \cite{peng2023generative}. We define the hierarchical Tucker ansatz, and then we describe the interpolation procedure that produces \(\tilde{P}\). The interpolation selects a set of pivots on each edge of the tree. The maximal-volume heuristics \cite{ryzhakov2024black} is one standard way for pivot selection, and our walkthrough of the first stage will assume that pivots are given for simplicity. For more detailed derivations, we refer readers to \cite{ballani2013black}.

\paragraph{Hierarchical Tucker ansatz}
The hierarchical Tucker format represents an order-\(d\) tensor by a binary tree of low-order components. In the complete binary tree case, a hierarchical Tucker tensor \(C \in \R^{n^{d}}\) has a low-rank condition for each node \(q = I_{k}^{(l)}\) with its parent node \(f = [d] \setminus I_{k}^{(l)}\). In particular, the unfolding along the bipartition \([d] = q \cup f\) has rank \(r_{f}\) and factors according to the equation
\begin{equation}\label{eqn: unfolding}
  C(i_{1}, \ldots, i_{d}) = \sum_{\gamma = 1}^{r_{f}} C_{q}(i_{q}, \gamma)\, C_{f}(\gamma, i_{f}).
\end{equation}
In diagram form, one writes
\begin{equation}\label{eqn: unfolding diagram}
  \HTFullBi{C} \;=\; \HTUnfold.
\end{equation}
When this holds at every node, \(C\) is represented by a tree of components, one bond per edge, in \(\mathcal{O}(d r^{3} + dnr)\) parameters rather than the \(n^{d}\) cost to store all entries. At a middle node \(q\), one can split \(C_{q}\) into tensor factorizations involving \(C_a\) and \(C_b\) for the two child nodes \(a\) and \(b\), which gives
\begin{equation}\label{eqn: htn forward map}
  C(i_{1}, \ldots, i_{d}) = \sum_{\alpha, \beta, \gamma} C_{a}(i_{a}, \alpha)\, C_{b}(i_{b}, \beta)\, F_{q}(\alpha, \beta, \gamma)\, C_{f}(\gamma, i_{f}),
\end{equation}
where \(F_{q}\) is the tensor component at \(q\), and \(C_{a}, C_{b}, C_{f}\) are the three subtree contractions. In diagram form, one has
\begin{equation}\label{eqn: htn forward map diagram}
  \HTFull{C} \;=\; \HTNodeMap{F_{q}}{C_{a}}{C_{b}}{C_{f}}\,.
\end{equation}

\paragraph{HT Interpolation}
In our setting, we seek \(\tilde{P}\) where \(\tilde{P}\) satisfies \Cref{eqn: htn forward map diagram} for every node \(q\) with \(\tilde{P}\) in place of \(C\).
The interpolation reconstructs the components of \(\tilde{P}\) from a small number of function queries by a CUR decomposition carried out edge by edge. An edge \(e = (q, v)\) belonging to the tree \(T\) splits \([d]\) into the variables on the side of \(q\) and those on the side of \(v\). We select \(r_{e}\) pivot assignments on each side. We write \(w^{(q \to v)}_{s}\) for the \(s\)-th pivot on the side of \(q\) and \(w^{(v \to q)}_{t}\) for the \(t\)-th pivot on the side of \(v\). We form a skeletonization of \(P\) on the \((q, v)\) edge as follows:
\begin{equation}\label{eqn: Z def}
  Z_{(q,v)}(s, t) = P\left(w^{(q \to v)}_{s},\, w^{(v \to q)}_{t}\right).
\end{equation}
Its singular value decomposition \(Z_{(q,v)} = U \Sigma W^{\top}\) fixes a gauge on the edge. We absorb \(\Sigma\) into the factor on the side of the root and leave the factor on the side of the leaves orthonormal. When \(v\) is the parent of \(q\), we set \(A_{q \to v} = U\) and \(A_{v \to q} = \Sigma W^{\top}\), and vice versa. 

With the edge gauges fixed, we form the right-hand side of the linear system for the desired component \(F_q\) by querying \(P\) at selected pivots. At a node \(q\), we abbreviate the three incoming edge factors by \(A_{a} := A_{a \to q}\), \(A_{b} := A_{b \to q}\), and \(A_{f} := A_{f \to q}\). At a middle node \(q\) with child nodes \(a, b\) and parent node \(f\), one has
\begin{equation}\label{eqn: B def}
  B_{q}(\alpha, \beta, \gamma) = P\left(w^{(a \to q)}_{\alpha},\, w^{(b \to q)}_{\beta},\, w^{(f \to q)}_{\gamma}\right),
\end{equation}
while at a leaf node one similarly has \(B_{q}(i, \gamma) = P\left(i,\, w^{(f \to q)}_{\gamma}\right)\). Then, one can obtain the tensor component \(F_{q}\) by solving the assembled linear system
\begin{equation}\label{eqn: stage one solve}
  \left( A_{a} \otimes A_{b} \otimes A_{f} \right) F_{q} = B_{q},
\end{equation}
where the three factors act on the corresponding bonds of \(F_{q}\). In diagram form, one has
\begin{equation}\label{eqn: stage one solve diagram}
  \HTSolveLHS \;=\; \HTSolveRHS\,.
\end{equation}
At the root, the parent factor is absent. At a leaf, the two child factors are absent, and \(A_{f}\, F_{q} = B_{q}\) determines the leaf component \(F_{q}(i, \gamma)\). We summarize the procedure in \Cref{alg: stage one}. The output is the signed reference \(\tilde{P}\) with components \((F_{q})_{q}\). Because the factorization interpolates \(P\) only at the chosen pivots, \(\tilde{P}\) can carry negative entries.

\begin{algorithm}[h]
  \caption{Hierarchical Tucker interpolation (stage one).}
  \label{alg: stage one}
  \begin{algorithmic}[1]
    \REQUIRE Function access to \(P \colon [n]^{d} \to \R\).
    \REQUIRE Complete binary tree of depth \(L\). Pivot sets and target ranks \(\{r_{e}\}\) on each edge \(e\).
    \FOR{each edge \(e = (q, v)\) of the tree}
      \STATE Form \(Z_{(q,v)}\) from \(P\) queries at the pivots of \(e\).
      \STATE Take the SVD \(Z_{(q,v)} = U \Sigma W^{\top}\). Absorb \(\Sigma\) into the factor on the side of the root, leaving the factor on the side of the leaves orthonormal.
    \ENDFOR
    \FOR{each node \(q\)}
      \STATE Form the pivot skeleton \(B_{q}\) from \(P\) queries.
      \STATE Solve \Cref{eqn: stage one solve} for the tensor component \(F_{q}\) (at the root, drop the parent factor, and at a leaf, drop both child factors and carry the physical axis).
    \ENDFOR
    \STATE \textbf{Output} the signed HT tensor \(\tilde{P}\) with components \((F_{q})_{q}\).
  \end{algorithmic}
\end{algorithm}

\subsection{Stage two: non-negative hierarchical Tucker fitting}\label{sec: NHT fitting}

The second stage fits a non-negative ansatz \(P_{G}\) to \(\tilde{P}\) under a variational formulation. The second stage finds a density through the signed measure while \(\tilde{P}\). As a pre-processing step, we apply a scalar scaling to each component of \(\tilde{P}\) to ensure that \(\lVert \tilde{P} \rVert_{F} = 1\) and that each component in \(\tilde{P}\) is the same in Frobenius norm before we start the second stage. The fitting is a generalization of the non-negative tensor-train procedure of \cite{tang2025variational} from the chain case to the binary tree case.

\paragraph{Loss function}
The objective combines the squared Frobenius data residual with a log barrier that keeps each component positive,
\begin{equation}\label{eqn: NHT loss}
  \ell(G) = \lVert P_{G} - \tilde{P} \rVert_{F}^{2} + \sum_{q} \mu_{q}\, \phi_{q}(G_{q}),
\end{equation}
where \(\mu_{q} > 0\) is the barrier coefficient at node \(q\), and \(\phi_{q}\) is the log barrier on the entries of the component,
\begin{equation}\label{eqn: barrier}
  \phi_{q}(G_{q}) =
  \begin{cases}
    -\sum_{\alpha, \beta, \gamma} \log G_{q}(\alpha, \beta, \gamma), & q \text{ a middle node},\\[2pt]
    -\sum_{i, \gamma} \log G_{q}(i, \gamma), & q \text{ a leaf},\\[2pt]
    -\sum_{\alpha, \beta} \log G_{q}(\alpha, \beta), & q \text{ the root}.
  \end{cases}
\end{equation}
The barrier replaces the positivity constraint with a smooth penalty. Minimizing the loss in \Cref{eqn: NHT loss} while adaptively decreasing the barrier coefficients \(\mu_{q}\) toward zero largely follows from the interior-point method \cite{nocedal1999numerical}. Writing \(\ell_{0}(G) := \lVert P_{G} - \tilde{P} \rVert_{F}^{2}\), one has
\begin{equation}\label{eqn: ell0 expand}
  \ell_{0}(G) = \langle P_{G}, P_{G} \rangle - 2 \langle P_{G}, \tilde{P} \rangle + \langle \tilde{P}, \tilde{P} \rangle.
\end{equation}

\paragraph{Alternating minimization}
The approach is summarized in \Cref{alg: fitting}. Essentially, one sweeps the components forward and then backward over the tree and anneals the barrier coefficients \(\mu_{q}\).
We minimize \(\ell\) one component at a time, holding the others fixed. One can see that \(\ell_{0}\) is a convex quadratic in \(G_{q}\) and the barrier \(\phi_{q}\) is strictly convex, and so \(\ell\) is strongly convex as a function of \(G_{q}\). At each visited component, we form the Newton direction \(\delta G_{q} = -(\nabla^{2}_{G_{q}} \ell)^{-1}(\nabla_{G_{q}} \ell)\) and update by \(G_{q} \gets G_{q} + \tau\, \delta G_{q}\). The step length \(\tau \in (0, 1]\) comes from a backtracking line search.

\paragraph{Gradient and Hessian by message passing}
The two inner products in \Cref{eqn: ell0 expand} that involve \(P_{G}\) are assembled by message passing over the tree. Fix a middle node \(q\) with child nodes \(a, b\) and parent \(f\). Let \(P_{G,a}(i_{a}, \alpha)\) be the contraction of the tensor components below child \(a\), and define \(P_{G,b}(i_{b}, \beta)\) likewise. Let \(P_{G,f}(i_{f}, \gamma)\) be the contraction of tensor components corresponding to non-descendants of \(q\). Let \(\tilde{P}_{a}, \tilde{P}_{b}, \tilde{P}_{f}\) be the matching subtree contractions of \(\tilde{P}\). The component \(G_{q}\) enters \(P_{G}\) through
\begin{equation}\label{eqn: component slice}
  P_{G}(i_{1}, \ldots, i_{d}) = \sum_{\alpha, \beta, \gamma} P_{G,a}(i_{a}, \alpha)\, P_{G,b}(i_{b}, \beta)\, G_{q}(\alpha, \beta, \gamma)\, P_{G,f}(i_{f}, \gamma),
\end{equation}
which is the forward map \Cref{eqn: htn forward map} with \(G_{q}\) in place of \(F_{q}\). In diagram, one has
\begin{equation}\label{eqn: component slice diagram}
  \HTFull{P_{G}} \;=\; \HTNodeMap{G_{q}}{P_{G,a}}{P_{G,b}}{P_{G,f}}\,.
\end{equation}
By message passing, one obtains the Gram matrices of these subtrees over their bonds,
\begin{equation}\label{eqn: self messages}
  M_{a}(\alpha, \alpha') = \sum_{i_{a}} P_{G,a}(i_{a}, \alpha)\, P_{G,a}(i_{a}, \alpha'),
\end{equation}
with \(M_{b}, M_{f}\) given by similar contractions over \(i_{b}, i_{f}\), respectively. Each is symmetric positive semidefinite of size \(\rho \times \rho\). Similarly, by message passing, one obtains
\begin{equation}\label{eqn: cross messages}
  L_{a}(\alpha, \alpha') = \sum_{i_{a}} P_{G,a}(i_{a}, \alpha)\, \tilde{P}_{a}(i_{a}, \alpha'),
\end{equation}
of size \(\rho \times r\), with \(L_{b}, L_{f}\) given by the same sums over \(i_{b}, i_{f}\). In diagram form, one has
\begin{equation}\label{eqn: messages diagram}
  M_{a} \;=\; \HTGram{P_{G,a}}{a}, \qquad L_{a} \;=\; \HTCross.
\end{equation}
At a leaf node, one has \(P_{G,a}(i_{a}, \alpha) = G_{a}(i_{a}, \alpha)\). Forming \(M_{\bullet}\) and \(L_{\bullet}\) is done recursively by message passing, i.e., by one sweep from the leaves to the root and back, and so the complexity for forming the terms is \(\mathcal{O}(d)\).

We write \((M_{a} \otimes M_{b} \otimes M_{f})\, G_{q}\) for the component obtained by applying \(M_{a}, M_{b}, M_{f}\) to the \(\alpha, \beta, \gamma\) axes of \(G_{q}\). We let \(\odot\) and \(\oslash\) denote the entrywise product and division. Substituting \Cref{eqn: component slice} into \Cref{eqn: ell0 expand} makes \(\ell_{0}\) quadratic in \(G_{q}\) through
\begin{equation}\label{eqn: inner products}
  \langle P_{G}, P_{G}\rangle = \langle G_{q},\, (M_{a} \otimes M_{b} \otimes M_{f})\, G_{q}\rangle, \quad
  \langle P_{G}, \tilde{P}\rangle = \langle G_{q},\, (L_{a} \otimes L_{b} \otimes L_{f})\, F_{q}\rangle,
\end{equation}
where \(\langle \cdot, \cdot\rangle\) is the entrywise inner product of two components. The gradient and Hessian at \(G_{q}\) then read
\begin{align}
  \nabla_{G_{q}} \ell &= 2\,(M_{a} \otimes M_{b} \otimes M_{f})\, G_{q} - 2\,(L_{a} \otimes L_{b} \otimes L_{f})\, F_{q} - \mu_{q}\,(1 \oslash G_{q}), \label{eqn: grad}\\
  \nabla^{2}_{G_{q}} \ell &= 2\,(M_{a} \otimes M_{b} \otimes M_{f}) + \mu_{q}\, \mathrm{diag}\!\left(1 \oslash (G_{q} \odot G_{q})\right). \label{eqn: hess}
\end{align}
The leaf and root components drop the absent bonds. At a leaf, only \(M_{f}\) acts, and at the root, \(M_{f}\) is absent. \Cref{sec: NHT fitting extra} exploits the leaf structure.

\begin{algorithm}[h]
  \caption{Non-negative hierarchical Tucker fitting (stage two).}
  \label{alg: fitting}
  \begin{algorithmic}[1]
    \REQUIRE Signed reference \(\tilde{P}\). Feasible non-negative initial \(P_{G}\) (\Cref{alg: warm start}).
    \REQUIRE Number of sweeps \(L_{\mathrm{it}}\). Barrier schedule \(\{\mu_{q}\}\).
    \FOR{\(\mathrm{it} = 1, \ldots, L_{\mathrm{it}}\)}
      \FOR{\(q\) forward over the tree, then backward}
        \STATE Refresh the messages \(M_{\bullet}, L_{\bullet}\) at \(q\). Update \(\mu_{q}\) by the barrier schedule.
        \STATE \(\delta G_{q} \gets -(\nabla^{2}_{\smash{G_{q}}} \ell)^{-1} (\nabla_{G_{q}} \ell)\) via \Cref{eqn: grad,eqn: hess}.
        \STATE \(\tau \gets\) backtracking line search on \(\ell\) along \(\delta G_{q}\).
        \STATE \(G_{q} \gets G_{q} + \tau\, \delta G_{q}\).
      \ENDFOR
    \ENDFOR
    \STATE \textbf{Output} the non-negative HT tensor \(P_{G}\).
  \end{algorithmic}
\end{algorithm}

\subsection{Acceleration of the Newton step}\label{sec: NHT fitting extra}

The dense Hessian of one component at a middle node acts on \(\rho^{3}\) entries. Inverting it costs \(\mathcal{O}(\rho^{9})\) and dominates the cost of \Cref{alg: fitting} at moderate rank. This subsection lowers the per-step cost and sets the barrier schedule.

\paragraph{Decoupling over the physical index}
At a leaf node \(q\), the loss \(\ell\) splits over the physical index \(i\). In other words, each slice \(G_{q}(i, :)\) couples only to itself in \(\ell\). Thus, computing the search direction in \(G_q\) splits into \(n\) independent \(\rho \times \rho\) systems. This costs \(\mathcal{O}(n \rho^{3})\) rather than the \(\mathcal{O}(n^{3}\rho^{3})\) of inverting the dense \((n\rho) \times (n\rho)\) Hessian.

\paragraph{Conjugate-gradient inner solve}
For middle and root nodes, we solve the Newton system by conjugate gradient (CG). We matricize the component as \(X(\gamma; (\alpha, \beta)) := G_{q}(\alpha, \beta, \gamma)\), of size \(\rho \times \rho^{2}\), with the parent bond as the row index. The action of the Hessian in \Cref{eqn: hess} on \(X\) is
\begin{equation}\label{eqn: matvec}
  X \;\longmapsto\; 2\, M_{f}\, X\, (M_{a} \otimes M_{b})^{\top} + D \odot X,
  \qquad D(\gamma; (\alpha, \beta)) = \frac{\mu_{q}}{G_{q}(\alpha, \beta, \gamma)^{2}}.
\end{equation}
The right multiplication is taken row by row: each row is reshaped into a \(\rho \times \rho\) matrix \(Y\) and mapped to \(M_{a}\, Y M_{b}^{\top}\), two matrix products. Each CG iteration then costs \(\mathcal{O}(\rho^{4})\), and our CG implementation truncates at \(300\) iterations for efficiency, and we allow early stopping when the linear system has been approximately solved. The case for \(q\) being a root node is the same in \Cref{eqn: matvec} if one omits the \(M_f\) term and the \(\gamma\) index.

\paragraph{Adaptive barrier and preconditioner}
We decrease the barrier coefficient in step with the data gradient. We apply a simple heuristic modified from \cite{nocedal2009adaptive}, with the following formula:
\begin{equation}\label{eqn: adaptive mu}
  \mu_{q}^{\mathrm{it}+1} = \min\!\left(\mu_{q}^{\mathrm{it}}, \; \tilde{\mu}_{q}^{\mathrm{it}}\right),
  \qquad
  \tilde{\mu}_{q}^{\mathrm{it}} = \frac{\sigma}{N_{q}} \sum_{\alpha, \beta, \gamma}  G_{q}(\alpha, \beta, \gamma) \cdot \left\lvert \nabla_{G_{q}} \ell_{0}(\alpha, \beta, \gamma)\right\rvert,
\end{equation}
with \(\sigma > 0\) a centering parameter (we use \(\sigma = 0.5\)) and \(N_{q}\) the entry count of \(G_{q}\). Notably, if one removes the absolute value in \Cref{eqn: adaptive mu}, one would exactly recover the heuristics in \cite{nocedal2009adaptive}. The schedule is non-increasing and floored at \(10^{-12}\). The CG preconditioner is the diagonal of the Hessian in \Cref{eqn: hess}: the sum of the data part \(\,(2\mathrm{diag}(M_{a})) \otimes (\mathrm{diag}(M_{b})) \otimes (\mathrm{diag}(M_{f}))\) and the barrier part \(\mathrm{diag}(\mu_{q}/G_{q}(\alpha, \beta, \gamma)^{2})\).

\subsection{Warm initialization}\label{sec: NHT warm start}

The Newton iteration of \Cref{alg: fitting} needs a feasible start. To allow for efficient optimization, we employ a warm initialization strategy. The reference itself would be the natural start, but its components carry signs, and simple strategies such as hard thresholding the negative entries would be a crude approximation. Instead, our proposal utilizes the gauge freedom in tensor networks to form a decoupled sweep that fits each non-negative component against the reference.

\paragraph{Gauge degree of freedom}
Let \(e\) be the edge between a node \(q\) and its parent \(f\), and let \((G_{q})_{q}\) be the components of an HT tensor on the tree. We gauge the edge by an orthogonal matrix \(Q_{e}\) of size \(\rho \times \rho\). We insert a pair of \(Q_e\) to apply to the parent bond of \(G_{q}\) and to the matching child bond of \(G_{f}\). The contraction over the edge meets both copies, and \(Q_{e}^{\top} Q_{e} = I\) cancels them. In tensor diagram, one writes
\begin{equation}\label{eqn: gauge cancel}
  \HTGaugeEdge \;=\; \HTGaugeEdgeRHS\,.
\end{equation}
Gauging every internal edge thus changes the components but not the underlying tensor. 
At node \(q\), the three gauges act as \((Q_{a} \otimes Q_{b} \otimes Q)\, G_{q}\), one per bond, with \(Q_{a}, Q_{b}\) on the child edges and \(Q\) on the parent edge.
In diagram form, the goal is to find \((G_q, Q, Q_a, Q_b)\) so that the following holds
\begin{equation}\label{eqn: warm local diagram}
  \HTWarmLHS \;\approx\; \HTWarmRHS\,.
\end{equation}

\paragraph{Main idea}

We propose a decoupled fitting strategy. Essentially, we decouple each individual component \(G_q\) from the environment \((G_v)_{v \not = q}\). Doing so allows for a faster bootstrap, as one is unlikely to find \(G_q\) stuck in a local minimum from a bad subspace as a result of the environment \((G_v)_{v \not = q}\). In simple terms, when optimizing \(G_q\), our strategy is to assume that the equation in \Cref{eqn: warm local diagram} holds \emph{exactly} for every node \(v\) other than \(q\). If that were true, one would have
\begin{equation}\label{eqn: leave-one-out}
    \begin{aligned}
        \ell_0(G) &= \lVert P_{G} - \tilde{P} \rVert_{F}^{2}\\
        &= \left\lVert (Q_{a} \otimes Q_{b} \otimes Q)\, G_{q} - F_{q} \right\rVert_{M}^{2},
    \end{aligned}
\end{equation}
where the \(M\)-norm is defined by \(\lVert X \rVert_{M}^{2} := \langle X,\, (M_{a} \otimes M_{b} \otimes M_{f})\, X \rangle\), where \(M_{a}, M_{b}, M_{f}\) are formed from the subtree contractions as in \Cref{eqn: self messages}, with the subtrees \(\tilde P_a, \tilde P_b, \tilde P_f\) in place of \(P_{G, a}\) in \Cref{eqn: self messages}. Thus, the main idea is to minimize the last line of \Cref{eqn: leave-one-out} over the \(G_q\) and \(Q\) variables. One can see that this amounts to a ``leave-one-out" type of bootstrapping technique, where \(G_q\) is solved by assuming that all other components have been adequately configured. While this is a heuristic procedure, we remark that the warm initialization strategy is paired with the main second stage algorithm in \Cref{alg: fitting} to ensure a good fitting.

\paragraph{Decoupled sweep}

We sweep nodes from the leaves to the root. At node \(q\), the child gauges \(Q_{a}, Q_{b}\) are already fixed, and we can only optimize the parent gauge variable \(Q\). With a local barrier coefficient \(\mu > 0\), we fit \(Q\) and a positive component \(G_{q}\) to the reference component,
\begin{equation}\label{eqn: warm local}
  \min_{G_{q} > 0, \; Q^{\top}Q = I} \; \left\lVert (Q_{a} \otimes Q_{b} \otimes Q)\, G_{q} - F_{q} \right\rVert_{M}^{2} - \mu \sum_{\alpha, \beta, \gamma} \log G_{q}(\alpha, \beta, \gamma).
\end{equation}
Similar to the stage two setting, the barrier keeps \(G_{q}\) entry-wise positive, and the barrier is adjusted dynamically. The formula is described by \Cref{eqn: adaptive mu} with \(\ell_0\) adapted to the quadratic term in \Cref{eqn: warm local}. We minimize over \(Q\) and \(G_q\) alternately, with at most forty rounds per node. Before the alternating minimization step, we initialize \(Q\) to be the identity matrix of size \(\rho \times \rho\), and we initialize \(G_{q}\) with i.i.d. entries from \(\mathrm{Unif}([0,1])\) with a simple scalar scaling to ensure \(\lVert G_{q} \rVert_{M} = \lVert F_{q} \rVert_{M}\).

\paragraph{\(Q\) step} With \(G_{q}\) fixed, only the loss term of \Cref{eqn: warm local} varies, and the problem reduces to a Procrustes problem weighted by \(M_{f}\). Write \(M_{a} = R_{a}^{\top} R_{a}\) for the Cholesky factorization, with \(R_{a}\) upper triangular, and likewise \(M_{b} = R_{b}^{\top} R_{b}\). Let \(\hat{C}_{q}\) denote \(G_{q}\) with \(Q_a, Q_b\) absorbed, followed by \(R_{a}\) and \(R_{b}\) absorbed on the child axes. Let \(\hat{F}_{q}\) denote \(F_{q}\) with \(R_{a}\) and \(R_{b}\) absorbed on the child axes. The weighting absorbs \(M_{a}\) and \(M_{b}\), so the child axes contract in the plain inner product. With \(X := \mathrm{mat}_{\gamma}(\hat{C}_{q})\) and \(Z := \mathrm{mat}_{\gamma}(\hat{F}_{q})\), where \(\mathrm{mat}_{\gamma}\) matricizes on the parent bond as in \Cref{sec: NHT fitting extra}, the loss term of \Cref{eqn: warm local} becomes
\begin{equation}\label{eqn: procrustes expand}
  \mathrm{tr}\left(\left(Q X - Z\right)^{\top} M_{f} \left(Q X - Z\right)\right)
  = \mathrm{tr}\left(X^{\top} Q^{\top} M_{f}\, Q X\right) - 2\, \mathrm{tr}\left(Q^{\top} K\right) + \mathrm{tr}\left(Z^{\top} M_{f}\, Z\right),
\end{equation}
with
\begin{equation}\label{eqn: procrustes K}
  K := M_{f}\, Z X^{\top}
    \;=\; \HTProcK\,.
\end{equation}
When \(M_{f}\) is the identity, one takes the singular value decomposition \(K = U \Sigma W^{\top}\), and \(Q = U W^{\top}\) would be the exact minimizer. We use \(Q = U W^{\top}\) as a heuristic, and we accept this update to \(Q\) only when the loss term decreases. At a leaf node, there are no child axes, so \(X = \mathrm{mat}_{\gamma}(G_{q})\) and \(Z = \mathrm{mat}_{\gamma}(F_{q})\). The root node has no parent edge and skips this step.

\paragraph{\(G\) step} With \(Q\) fixed, \Cref{eqn: warm local} over \(G_{q}\) is the node problem of \Cref{sec: NHT fitting} in a rotated frame: the Hessian is \Cref{eqn: hess} with each message conjugated by the gauge on its edge, so \(M_{a}\) becomes \(Q_{a}^{\top} M_{a} Q_{a}\), and likewise for \(M_{b}\) and \(M_{f}\). The conjugation preserves the Kronecker-plus-diagonal structure, so the matrix-free CG of \Cref{sec: NHT fitting extra} applies, similarly with the Hessian diagonal as the preconditioner. We take one interior-point Newton step with a backtracking line search on the objective of \Cref{eqn: warm local}.

\begin{algorithm}[h]
  \caption{Gauge warm initialization.}
  \label{alg: warm start}
  \begin{algorithmic}[1]
    \REQUIRE Signed reference \(\tilde{P}\) with components \((F_{q})_{q}\).  Barrier parameter \(\nu\) and floor \(\mu_{\min}\).
    \STATE Compute messages from the original signed components.
    \FOR{\(q\) from leaves to root}
      \STATE Solve \Cref{eqn: warm local} by alternating the \(Q\) step and the \(C\) step, each accepted only if the loss term of \Cref{eqn: warm local} decreases. Record \(Q\) on the parent edge and the component \(G_{q}\).
    \ENDFOR
    \STATE \textbf{Output} the non-negative HT factors \((G_{q})_{q}\).
  \end{algorithmic}
\end{algorithm}

\paragraph{Algorithmic summary} We summarize our approach in \Cref{alg: warm start}. By alternatively optimizing \(G_q\) and \(Q\), we can obtain good initial approximations in practice. 

\section{Numerical experiments}\label{sec: numerics}

We test the proposed NHT compression on both access models of \Cref{sec: alg}. \Cref{sec: VI numerics} treats the variational inference case, where one queries entries of \(P\). \Cref{sec: DE numerics} treats the density estimation case, where one is given samples of \(P\). The experiments evaluate the second stage. Every method starts from the same gauge warm initialization of \Cref{alg: warm start}, followed by one further scalar rescale to the loss-minimizing scale, and fits the same signed reference \(\tilde{P}\). Each panel therefore isolates the fitting iteration: the multiplicative-update benchmark of \Cref{alg: lee-seung} against three variants of the Newton iteration of \Cref{alg: fitting}. The variants differ in the barrier schedule (fixed or adaptive) and the CG preconditioning of \Cref{sec: NHT fitting extra}.

Lastly, to justify the choice of the warm initialization procedure in \Cref{alg: warm start} and the choice of the NHT ansatz, we conduct an ablation study in \Cref{sec: ablation}. The result shows that the choice of NHT leads to better fitting error than the NTT ansatz in the 2D lattice model considered in \Cref{sec: VI numerics}. Moreover, we show that \Cref{alg: warm start} is more suitable for NHT fitting than using the multiplicative update algorithm in \Cref{alg: lee-seung} for initialization.

\paragraph{Protocol}
Each panel plots the relative squared Frobenius loss \(\lVert P_{G} - \tilde{P}\rVert_{F}^{2}/\lVert \tilde{P}\rVert_{F}^{2}\) against wall-clock time from the shared start. Each method receives a \(100\)-second fitting budget after the shared warm initialization. The adaptive schedule starts each node from the warm \(\mu_{\mathrm{ref}}\), and the fixed schedule starts at \(\mu = 5 \times 10^{-4}\) and halves each sweep, down to the same \(10^{-12}\) floor. The benchmark runs with floor \(\varepsilon = 10^{-10}\) under the same wall-clock budget. All computations are done with the CPU of an M5 Pro chip on a MacBook Pro.

\paragraph{Benchmark} For benchmark, we choose the multiplicative update algorithm \cite{lee2000algorithms}, and we remark that the implementation is adapted to the tree case from the tensor-train version of \cite{shcherbakova2019nonnegative}. Each component is rescaled entrywise by the ratio of the cross-gradient \(\nabla_{G_{q}} \langle P_{G}, \tilde{P}\rangle\) to half the self-gradient \(\nabla_{G_{q}} \langle P_{G}, P_{G}\rangle\), each floored at \(\varepsilon > 0\). The floor keeps the update non-negative even where the signed reference makes the cross-gradient negative. The benchmark replaces the Newton step of \Cref{alg: fitting} with this update. Because the update can be non-monotone against the signed reference, we reject any sweep that raises \(\ell_{0}\). \Cref{alg: lee-seung} summarizes the approach.

\begin{algorithm}[h]
  \caption{Multiplicative update benchmark.}
  \label{alg: lee-seung}
  \begin{algorithmic}[1]
    \REQUIRE Signed reference \(\tilde{P}\). Feasible non-negative \(P_{G}\). Number of sweeps \(L_{\mathrm{it}}\). Floor \(\varepsilon\).
    \FOR{\(\mathrm{it} = 1, \ldots, L_{\mathrm{it}}\)}
      \FOR{\(q\) forward over the tree, then backward}
        \STATE \(G_{q} \gets G_{q} \odot \dfrac{\max(\nabla_{G_{q}} \langle P_{G}, \tilde{P}\rangle, \; \varepsilon)}{\max(\tfrac{1}{2}\nabla_{G_{q}} \langle P_{G}, P_{G}\rangle, \; \varepsilon)}\) \quad (entrywise).
      \ENDFOR
      \STATE If \(\ell_{0}\) increased over the sweep, restore the components from before the sweep and stop.
    \ENDFOR
    \STATE \textbf{Output} the non-negative HT factors \((G_{q})_{q}\).
  \end{algorithmic}
\end{algorithm}

\subsection{Variational inference}\label{sec: VI numerics}

In the variational inference case, the first stage is the interpolation of \Cref{alg: stage one}, for which we take the approach in \cite{ryzhakov2024black}. To measure the accuracy of the fit, we draw \(10^{5}\) held-out multi-indices uniformly from \([n]^{d}\). We report the relative \(\ell^{2}\) error \(\lVert P_{G} - P\rVert / \lVert P\rVert\) over those indices.

\begin{figure}[tbp]
  \centering
  \subcaptionbox{\label{fig: vi gl}}[0.86\linewidth]{
    \includegraphics[width=\linewidth]{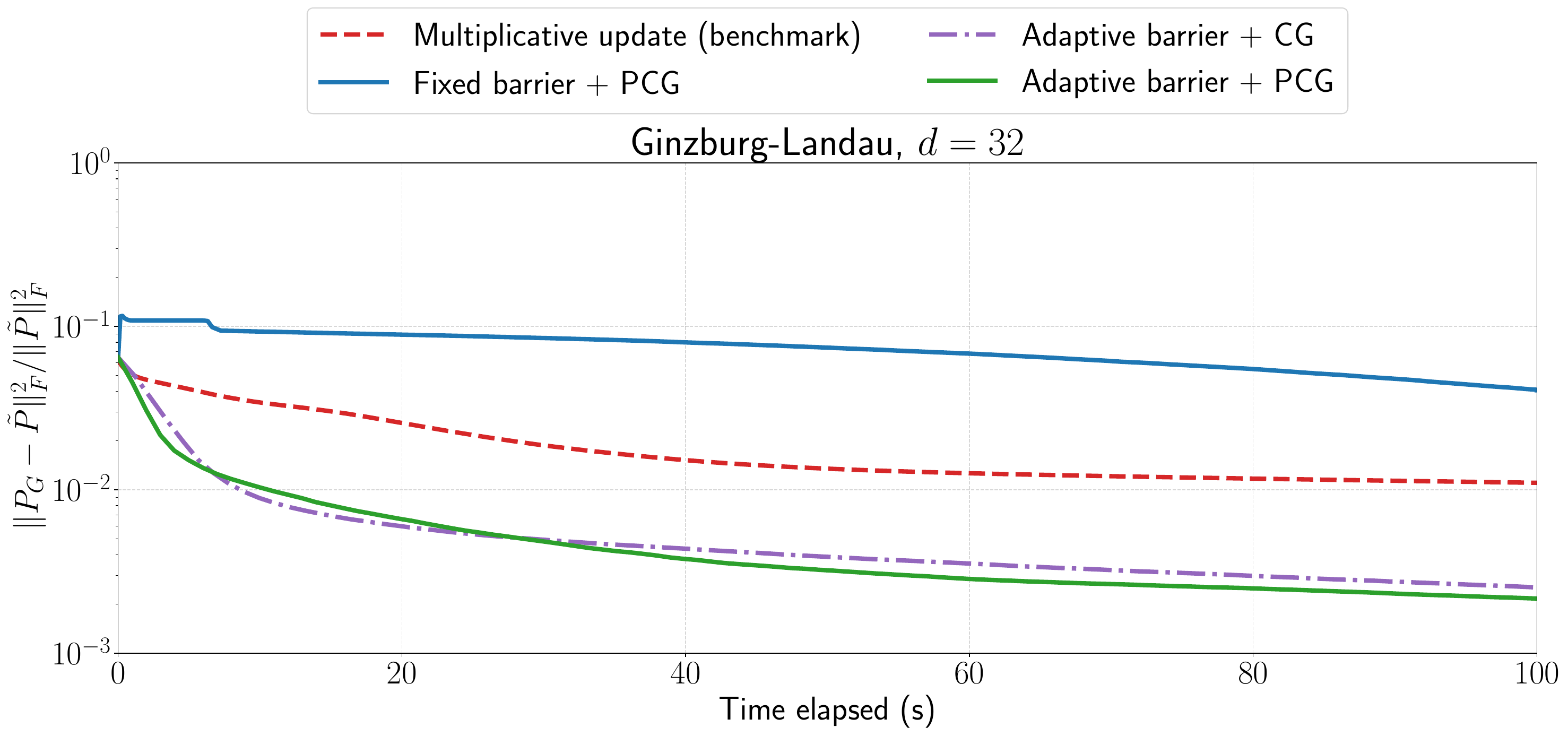}}\\[6pt]
  \subcaptionbox{\label{fig: vi cauchy}}[0.86\linewidth]{
    \includegraphics[width=\linewidth]{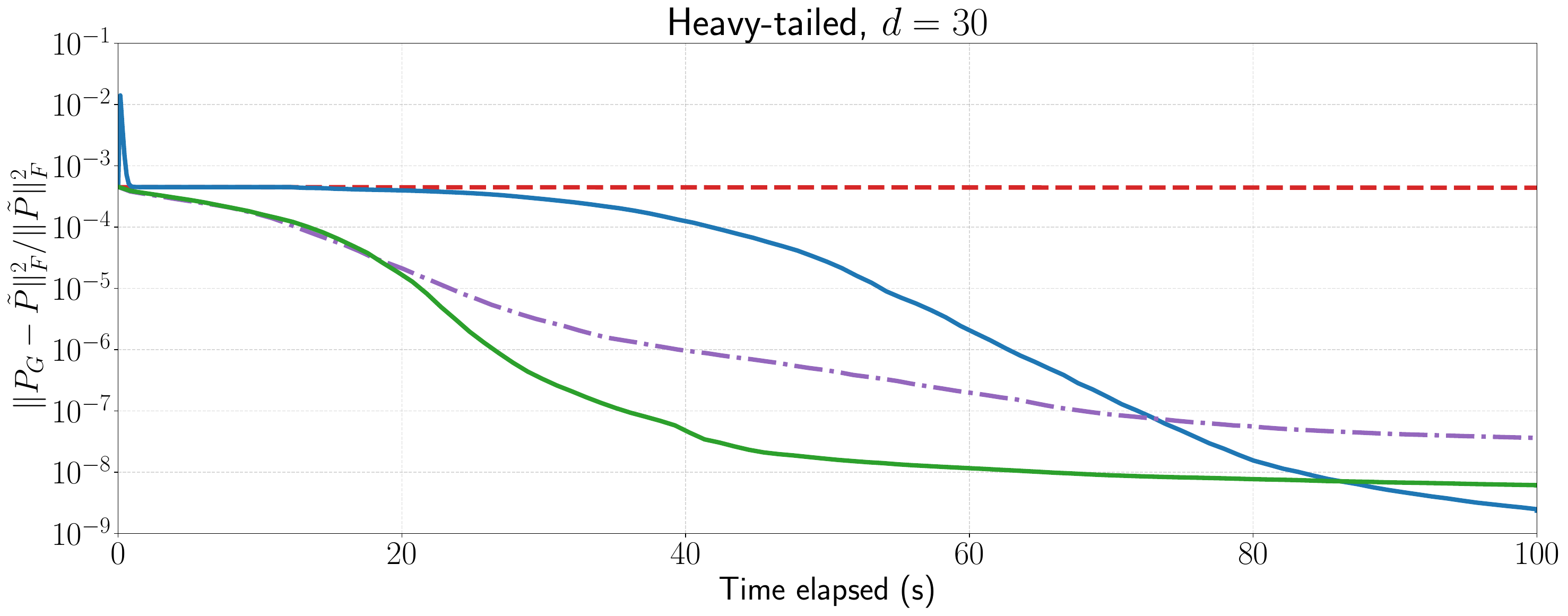}}
  \caption{NHT fitting in the variational inference case. Each panel plots the relative squared Frobenius loss against wall-clock time, with every curve started from the same gauge warm initialization. We compare the Newton iteration against the multiplicative-update benchmark. The three Newton variants vary the barrier schedule and the preconditioning of the CG inner solve. The adaptive-barrier preconditioned iteration substantially improves on the multiplicative-update benchmark in both examples.}
  \label{fig: vi numerics}
\end{figure}

\paragraph{Example 1: Ginzburg-Landau model}
We discretize the model of \Cref{sec: alg} on a \(4 \times 8\) lattice, so \(d = 32\), and we place \(n = 21\) grid points \((x_{i})_{i=1}^{n}\) on \([-2, 2]\). The distribution tensor is
\begin{equation}\label{eqn: GL formula}
  P(i_{1}, \ldots, i_{d}) = \exp\!\left(-\beta\left(\frac{\gamma}{2d}\sum_{v \sim w} (x_{i_{v}} - x_{i_{w}})^{2} + \frac{\lambda}{d}\sum_{v} (x_{i_{v}}^{2} - 1)^{2}\right)\right),
\end{equation}
where \(v \sim w\) ranges over the nearest-neighbor pairs of the periodic lattice. Here, \(\beta\) is the inverse temperature, \(\gamma\) controls the correlation strength, and \(\lambda\) controls the double-well strength. We take \(\beta = 16\), \(\gamma = 0.10\), and \(\lambda = 0.25\). In the first stage, the interpolation with a maximal internal rank of \(14\) compresses \(P\) into \(\tilde{P}\) with a relative error of \(7.8 \times 10^{-2}\) on the held-out indices. In the second stage, we fit an NHT ansatz \(P_{G}\) of internal rank \(14\). \Cref{fig: vi gl} shows that the adaptive-barrier preconditioned Newton iteration reaches a relative squared Frobenius loss of \(2.2 \times 10^{-3}\), which is below the multiplicative-update loss at the same wall-clock time. On the held-out indices, the relative error of \(P_{G}\) against \(P\) is \(9.5 \times 10^{-2}\) for the preconditioned case.

\paragraph{Example 2: heavy-tailed model}
We take \(n = 50\) grid points \((x_{i})_{i=1}^{n}\) on \([0, 2]\) and the heavy-tailed distribution
\begin{equation}\label{eqn: cauchy numerics}
  P(i_{1}, \ldots, i_{d}) = \frac{1}{1 + x_{i_{1}}^{2} + \cdots + x_{i_{d}}^{2}}
\end{equation}
where \(d = 30\). This is a multivariate analogue of the Cauchy distribution and is symmetric in all variables. In the first stage, the interpolation with a maximal internal rank of \(10\) compresses \(P\) into \(\tilde{P}\) with a relative error of \(3.2 \times 10^{-6}\) on the held-out indices. In the second stage, we fit an NHT ansatz of internal rank \(10\). \Cref{fig: vi cauchy} shows that the adaptive-barrier preconditioned Newton iteration reaches a relative squared Frobenius loss of \(6.1 \times 10^{-9}\), close to the lowest of the fitting variants. On the held-out indices, the relative error of \(P_{G}\) against \(P\) is \(7.2 \times 10^{-5}\).

\subsection{Density estimation}\label{sec: DE numerics}

In the density estimation case, the first stage replaces the function queries of \Cref{sec: HT} by empirical means over the samples \cite{peng2023generative}. Since \(\tilde{P}\) converges to \(P\) only at the Monte Carlo rate, it carries non-negligible negative entries, and so the fitting loss from a non-negative ansatz will not reach zero. We therefore also report the negative log-likelihood (NLL),
\begin{equation}\label{eqn: nll}
  \mathrm{NLL}(P_{G}) = -\frac{1}{N} \sum_{j=1}^{N} \log\frac{P_{G}(y^{(j)})}{Z_{G}}, \qquad Z_{G} = \sum_{i_{1}, \ldots, i_{d}} P_{G}(i_{1}, \ldots, i_{d}),
\end{equation}
evaluated on the samples \((y^{(j)})_{j=1}^{N}\). The normalizing constant \(Z_{G}\) is computed in \(\mathcal{O}(d)\) operations by summing each leaf component over its physical index and contracting the tree. The fitting succeeds when \(\mathrm{NLL}(P_{G})\) is close to the NLL of the ground-truth \(P\).

\begin{figure}[tbp]
  \centering
  \subcaptionbox{\label{fig: de ising}}[0.86\linewidth]{
    \includegraphics[width=\linewidth]{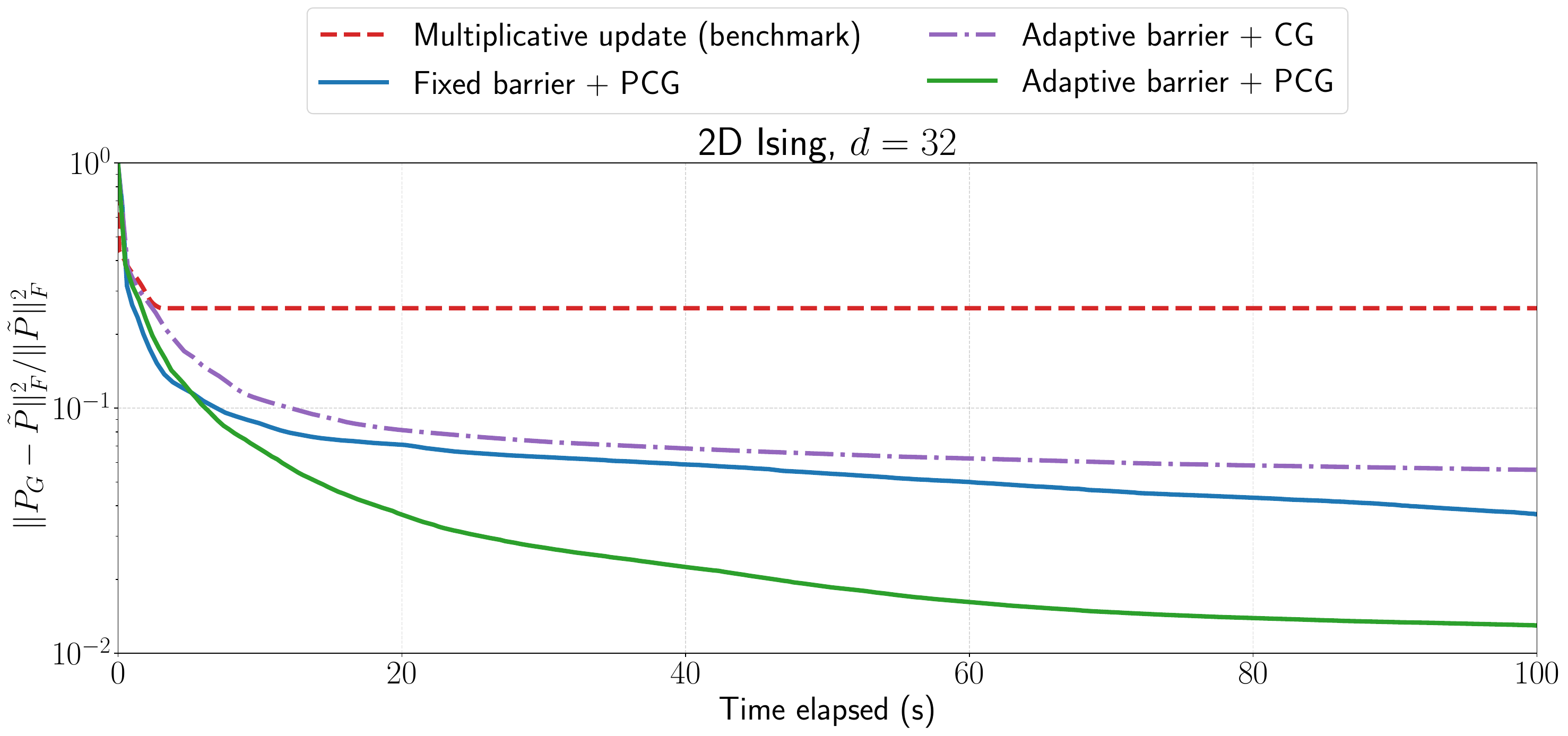}}\\[6pt]
  \subcaptionbox{\label{fig: de tfim}}[0.86\linewidth]{
    \includegraphics[width=\linewidth]{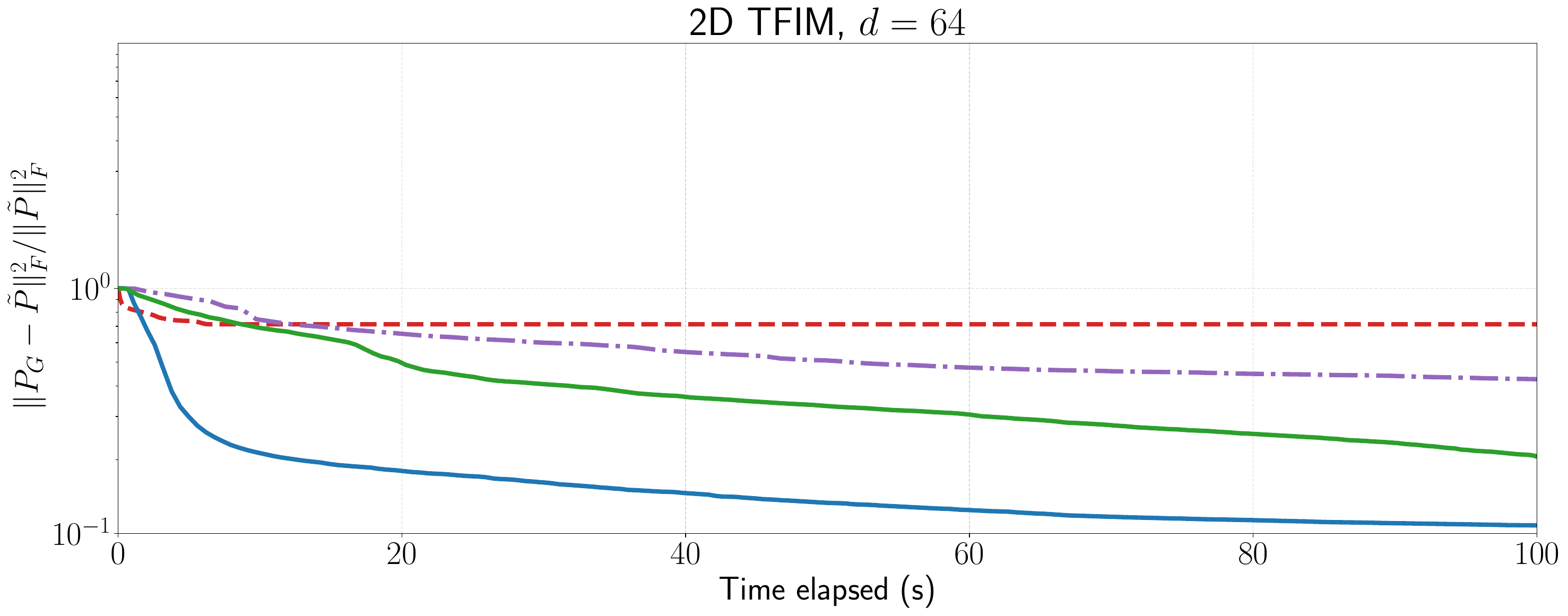}}
  \caption{NHT fitting in the density estimation case. As in \Cref{fig: vi numerics}, every curve, including the multiplicative-update benchmark, starts from the same gauge warm initialization. Each panel plots the relative squared Frobenius loss against wall-clock time. The three Newton variants vary the barrier schedule and the preconditioning of the CG inner solve. The Newton variants substantially improve on the multiplicative-update benchmark, with the adaptive-barrier preconditioned iteration strongest on the Ising example.}
  \label{fig: de numerics}
\end{figure}

\paragraph{Example 3: two-dimensional Ising model}
We consider a ferromagnetic Ising model on a \(4 \times 8\) torus, so \(d = 32\) with \(n = 2\). The distribution is
\begin{equation}\label{eqn: ising numerics}
  P(i_{1}, \ldots, i_{d}) \propto \exp\!\left(\beta \sum_{v \sim w} s_{i_{v}} s_{i_{w}}\right), \quad s_{i} \in \{-1, +1\}, \quad \beta = 0.10,
\end{equation}
with \(v \sim w\) ranging over the nearest-neighbor pairs of the torus. We draw \(N = 5 \times 10^{5}\) samples by the Wolff cluster algorithm \cite{liu2001monte}. We fit an NHT ansatz of internal rank \(16\). \Cref{fig: de ising} shows that the adaptive-barrier preconditioned Newton iteration has a relative squared Frobenius loss of \(1.3 \times 10^{-2}\). We test accuracy with the NLL. The exact NLL of \(P\) is \(21.85\). The NLL of \(P_{G}\) on the samples is \(21.89\).

\paragraph{Example 4: two-dimensional transverse-field Ising model}
We consider the ground state of a two-dimensional transverse-field Ising model on an \(8 \times 8\) torus with \(d = 64\). The Hamiltonian is \(H = -\sum_{v \sim w} Z_{v} Z_{w} - h \sum_{v} X_{v}\), where \(Z_{v}, X_{v}\) are the Pauli operators on site \(v\). We take \(h = 4\). The distribution tensor is \(P(i_{1}, \ldots, i_{d}) = \lvert \ket{\psi}(i_{1}, \ldots, i_{d}) \rvert^{2}\), the Born distribution of the ground state \(\ket{\psi}\). We obtain \(\ket{\psi}\) by the density-matrix renormalization group \cite{white1992density} at maximal bond dimension \(150\), and we form the first-stage sketch from \(8 \times 10^{6}\) samples from the Born distribution by exact sampling of the matrix product state. We fit an NHT ansatz of internal rank \(16\). \Cref{fig: de tfim} shows that the adaptive-barrier preconditioned iteration improves over the multiplicative-update benchmark and reaches a loss of \(2.1 \times 10^{-1}\), while the fixed-barrier preconditioned variant reaches the lowest loss on this example. On \(2 \times 10^{5}\) diagnostic samples, the NLL of \(P\) is \(43.02\), and the NLL of \(P_{G}\) is \(43.52\).

\subsection{Ablation studies}\label{sec: ablation}

This subsection aims to address two natural architectural questions regarding this work. First is about the use of the hierarchical Tucker ansatz in NHT over the NTT ansatz. Second is about the use of the warm initialization algorithm in \Cref{alg: warm start} over the more standard multiplicative update algorithm in \Cref{alg: lee-seung}. To address the first question, we shall compare the NHT ansatz with the NTT ansatz in the 2D Ginzburg-Landau model considered in \Cref{sec: VI numerics}. For simplicity, in addressing the second question, we also fold in the warm initialization strategy comparison in that example. For fair comparison, all methods considered in this subsection use an adaptive barrier schedule along with the preconditioned CG algorithm during the Newton step.

\begin{figure}[tbp]
  \centering
    \includegraphics[width=0.86\linewidth]{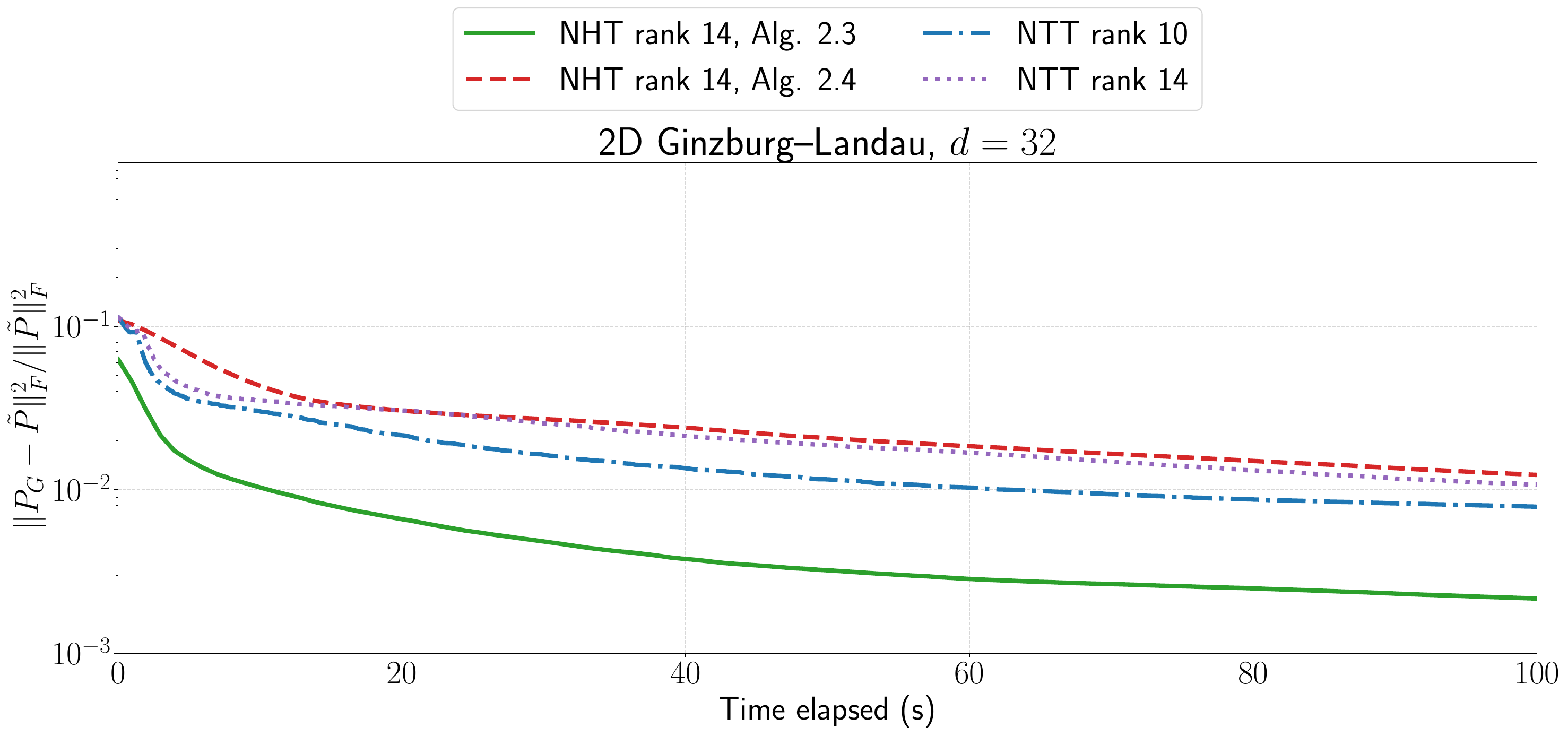}
  \caption{Ablation study that compares the default NHT ansatz implementation against other architectural alternatives. The experiment is on the 2D Ginzburg-Landau setting considered in \Cref{sec: VI numerics}. The result shows that NHT is better at compressing the ground truth density than NTT. Moreover, the result shows that the warm initialization procedure in \Cref{alg: warm start} is better for downstream NHT fitting than \Cref{alg: lee-seung}.}
  \label{fig: ablation}
\end{figure}

The experiment setup is simple. The variational inference problem is as in \Cref{eqn: GL formula}, where dimension is \(d=32\), and the problem is a 2D \(4 \times 8\) lattice. In the first stage, we use CUR-type interpolation algorithms to compress \(P\) into a signed TT target \(\tilde{P}_{TT}\) and a signed HT target \(\tilde{P}_{HT}\). In the second stage for \(\tilde{P}_{HT}\), we use the NHT ansatz to fit \(\tilde{P}_{HT}\), and we include the result of using \Cref{alg: warm start} and \Cref{alg: lee-seung} as the warm initialization method. In the second stage for \(\tilde{P}_{TT}\), we follow the procedure in \cite{tang2025variational} to use a NTT ansatz to fit \(\tilde{P}_{TT}\). In this setting, we use rank \(\rho_{HT} = 14\) for the NHT ansatz and we use rank \(\rho_{TT} \in \{10, 14\}\) for the NTT ansatz. By direct calculation, one can see that the \(\rho_{TT} = 10\) case gives the NTT ansatz roughly the same number of parameters as the NHT ansatz with \(\rho_{HT} = 14\). Moreover, the NTT ansatz under \(\rho_{TT} = 14\) has roughly double the number of parameters as the NHT ansatz with \(\rho_{HT} = 14\). We show the result in \Cref{fig: ablation} for the squared relative error:
\[
\ell_{0}(G) = \frac{\lVert P_G - \tilde{P} \rVert^2 }{\lVert\tilde{P} \rVert^2},
\]
where we take \(\tilde{P}\) to be \(\tilde{P}_{TT}\) for the NTT ansatz and we take \(\tilde{P}\) to be \(\tilde{P}_{HT}\) for the NHT ansatz. Since the dimension is still relatively small at \(d =32\), we have \(\tilde{P}_{TT} \approx \tilde{P}_{HT}\) as \(d\)-tensors. We remark that the first stage output \(\tilde{P}_{TT}\) and \(\tilde{P}_{HT}\) are both fairly accurate, and so one can simply view \(\lVert P_G - \tilde{P} \rVert^2 /\lVert\tilde{P} \rVert^2\) as a high-fidelity approximation of \(\lVert P_G - P \rVert^2 /\lVert P \rVert^2\), where \(P\) is the ground truth probability tensor. On the held-out indices, the relative error of \(P_{G}\) against \(P\) is \(9.5 \times 10^{-2}\) for NHT under \Cref{alg: warm start}. For NHT under \Cref{alg: lee-seung}, NTT with \(\rho_{TT} = 10\), and NTT with \(\rho_{TT} = 14\),  we respectively obtain a relative error of \(1.5\times 10^{-1}\), \(1.8\times 10^{-1}\) and \(1.6\times 10^{-1}\). 

From \Cref{fig: ablation}, one sees that the NTT ansatz is less suitable for compressing the 2D Ginzburg-Landau model than NHT, and the gap does not close even when one gives NTT double the number of parameters as the NHT model. Similarly, the result in \Cref{fig: ablation} shows a clear separation between the performance of \Cref{alg: warm start} and \Cref{alg: lee-seung} as warm initialization strategies. While we only demonstrate the performance for the 2D G-L case, we remark that the choice of the NHT ansatz along with \Cref{alg: warm start} as warm initialization consistently leads to the best fitting among all architectural designs we have considered. 

\section{Discussion}\label{sec: conclusion}

We introduced an end-to-end approach for compressing a high-dimensional distribution with the non-negative hierarchical tensor ansatz. The approach enjoys fast convergence
and complements the existing non-negative tensor-train algorithm \cite{tang2025variational} for target densities that do not satisfy an area law, notably from 2D lattice models. An interesting future direction is to use the NHT ansatz for the coefficient tensor in the functional hierarchical tensor ansatz \cite{tang2024solving}. Future work can consider combining the NHT fitting task with additional requirements, such as moment conservation.

\appendix
\section{Extension to general tree tensor networks}\label{sec: appendix TTNS}

\Cref{sec: alg} is written for the complete binary tree. This appendix extends the two-stage method to general tree tensor networks \cite{tang2023generative}, and we mainly follow the exposition in \cite{tang2025wavelet}.

\paragraph{Tree structure notation}
A tree graph \(T = (V, E)\) is a connected undirected graph without cycles. For a node \(q \in V\), \(\mathcal{N}(q)\) is the set of its neighbors, and \(\mathcal{E}(q)\) is the set of its incident edges. Removing an edge \(e = (v, q)\) from \(E\) leaves two connected components, and \(v \to q\) denotes the side that contains \(v\). One node of \(T\) is prescribed as the root, and the root orients the tree: every other node has a parent, the neighbor on the root side, and its remaining neighbors are its children.

In simple terms, a tree tensor network stores a \(d\)-tensor by one low-order component per node of a tree, contracted along the edges. A node that carries a variable is external, and a node that carries none is internal.

\begin{definition}[Tree tensor network]\label{def: ttns}
Let \(T = (V, E)\) be a tree graph with ranks \(\{\rho_{e}\}_{e \in E}\), and let \(V_{\mathrm{ext}} \subseteq V\) be the set of external nodes, with \(d = \lvert V_{\mathrm{ext}} \rvert\) and all other nodes internal. We label the nodes so that \(V_{\mathrm{ext}} = [d]\). One tensor component sits on each node, of shape
\begin{equation}\label{eqn: ttns shapes}
  G_{q} \colon [n] \times \prod_{e \in \mathcal{E}(q)} [\rho_{e}] \to \R \;\; \text{at an external node},
  \qquad
  G_{q} \colon \prod_{e \in \mathcal{E}(q)} [\rho_{e}] \to \R \;\; \text{at an internal node}.
\end{equation}
A \(d\)-tensor \(P_{G}\) is a tree tensor network over \(T\) with components \((G_{q})_{q \in V}\) when
\begin{equation}\label{eqn: ttns contraction}
  P_{G}(i_{1}, \ldots, i_{d}) = \sum_{\alpha_{E}}\, \prod_{q \in V_{\mathrm{ext}}} G_{q}\left(i_{q}, \alpha_{\mathcal{E}(q)}\right) \prod_{q \in V \setminus V_{\mathrm{ext}}} G_{q}\left(\alpha_{\mathcal{E}(q)}\right).
\end{equation}
\end{definition}

\begin{figure}[h]
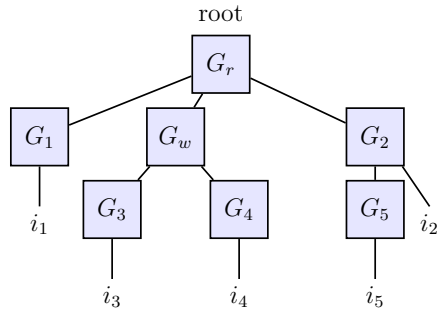

  \centering
  \scalebox{0.8}{\HTTTNSFigure}
  \caption{A tree tensor network on seven nodes. The external nodes carry the physical legs \(i_{1}, \ldots, i_{5}\), the internal nodes \(r\) and \(w\) carry none, and \(r\) is the prescribed root. The tensor component \(G_2\) carries both bonds and a physical leg, which has no counterpart in the complete binary tree setting in the main text.}
  \label{fig: ttns}
\end{figure}

Here \(i_{q} \in [n]\) is the physical index of an external node \(q\), \(\alpha_{e}\) is the bond on edge \(e\), \(\alpha_{\mathcal{E}(q)}\) collects the bonds on the edges of \(q\), and \(\alpha_{E}\) collects the bonds on all edges. We also write \(i_{v \to q}\) for the physical subvector on the external nodes of the side \(v \to q\). As in \Cref{sec: alg}, we write the ranks as a uniform \(\rho\). The signed reference \(\tilde{P}\) has components \((F_{q})_{q}\) of the same shapes, with rank \(r\) in place of \(\rho\) on every bond. \Cref{fig: ttns} fixes the node convention: an external node carries a physical leg, and an internal node does not. In the diagram equations below, a thick open leg is the grouped bond \(\alpha_{\mathcal{E}(q)}\), all edges of \(q\) drawn as one line, while a contraction is always a thin line. \Cref{sec: alg} is the special case in which the external nodes are the leaves and the three directions into a middle node are \(a \to q\), \(b \to q\), \(f \to q\).

\paragraph{Contractions at a node}
Both stages solve for one component at a time, and the equations at a node \(q\) involve the rest of the tree through one contraction per neighbor. For each neighbor \(v \in \mathcal{N}(q)\), we construct \(P_{G, v \to q}\) by contracting all tensor components \(G_{w}\) for \(w \in v \to q\). The result carries the physical subvector of its side and the bond of its edge, \(P_{G, v \to q}\left(i_{v \to q}, \alpha_{(v,q)}\right)\). The component \(G_{q}\) enters \(P_{G}\) only through these contractions: at an external node,
\begin{equation}\label{eqn: ttns node identity}
  P_{G}(i_{1}, \ldots, i_{d}) = \sum_{\alpha_{\mathcal{E}(q)}} G_{q}\left(i_{q}, \alpha_{\mathcal{E}(q)}\right) \prod_{v \in \mathcal{N}(q)} P_{G, v \to q}\left(i_{v \to q}, \alpha_{(v,q)}\right),
\end{equation}
and at an internal node, the same identity holds without the \(i_{q}\) index. \Cref{eqn: ttns node identity} is the analogue of \Cref{eqn: component slice}. The reference \(\tilde{P}\) has the contractions \(\tilde{P}_{v \to q}\), defined the same way from \((F_{q})_{q}\).

\paragraph{Stage one}
The interpolation selects pivots edge by edge and then solves node by node. On each edge \(e = (q, v)\) we select \(r_{e}\) pivot assignments per side, written \(w^{(q \to v)}_{s}\) and \(w^{(v \to q)}_{t}\), and we form the skeleton \(Z_{(q,v)}(s, t) = P(w^{(q \to v)}_{s}, w^{(v \to q)}_{t})\) from queries of \(P\). Its singular value decomposition \(Z_{(q,v)} = U \Sigma W^{\top}\) fixes the gauge on the edge, and the factor on the root side absorbs \(\Sigma\): when \(v\) is the parent of \(q\), \(A_{q \to v} = U\) and \(A_{v \to q} = \Sigma W^{\top}\), and when \(q\) is the parent of \(v\), \(A_{q \to v} = U \Sigma\) and \(A_{v \to q} = W^{\top}\). The pivot skeleton at a node evaluates \(P\) at one pivot per incident edge, with the physical index free,
\begin{equation}\label{eqn: ttns B def}
  B_{q}\left(i_{q}, \alpha_{\mathcal{E}(q)}\right) = P\Bigl(i_{q},\, \bigl(w^{(v \to q)}_{\alpha_{(v,q)}}\bigr)_{v \in \mathcal{N}(q)}\Bigr)
\end{equation}
at an external node, and without the \(i_{q}\) index at an internal node. Evaluating \Cref{eqn: ttns node identity} for \(\tilde{P}\) at these pivots replaces each contraction \(\tilde{P}_{v \to q}\) by the edge factor \(A_{v \to q}\). With \(A_{q} := \bigotimes_{v \in \mathcal{N}(q)} A_{v \to q}\), the tensor component therefore solves the linear system
\begin{equation}\label{eqn: ttns solve}
  A_{q}\, F_{q} = B_{q}.
\end{equation}
In diagram form, at an external node,
\begin{equation}\label{eqn: ttns solve diagram}
  \HTTTNSSolveLHS \;=\; \HTTTNSSolveRHS\,.
\end{equation}
At an internal node, the physical leg is absent. The root needs no special case, since \(\mathcal{E}(q)\) ranges only over the incident edges. \Cref{alg: stage one} runs with these replacements, and its output is the signed reference \(\tilde{P}\) with components \((F_{q})_{q}\).

\paragraph{Stage two}
Each node problem of the fitting stage is a convex quadratic plus a barrier, assembled from one message per incident edge. The loss keeps the form of \Cref{eqn: NHT loss},
\begin{equation}\label{eqn: ttns loss}
  \ell(G) = \lVert P_{G} - \tilde{P} \rVert_{F}^{2} + \sum_{q \in V} \mu_{q}\, \phi_{q}(G_{q}),
\end{equation}
with \(\phi_{q}\) the log barrier on all entries of \(G_{q}\), as in \Cref{eqn: barrier}. Each incident edge carries a self-message and a cross-message, built from the contractions \(P_{G, v \to q}\) and \(\tilde{P}_{v \to q}\),
\begin{align}
  M_{v \to q}(\alpha, \alpha') &= \sum_{i_{v \to q}} P_{G, v \to q}\left(i_{v \to q}, \alpha\right) P_{G, v \to q}\left(i_{v \to q}, \alpha'\right), \label{eqn: ttns messages}\\
  L_{v \to q}(\alpha, \alpha') &= \sum_{i_{v \to q}} P_{G, v \to q}\left(i_{v \to q}, \alpha\right) \tilde{P}_{v \to q}\left(i_{v \to q}, \alpha'\right), \label{eqn: ttns cross}
\end{align}
of sizes \(\rho \times \rho\) and \(\rho \times r\), the analogues of \Cref{eqn: self messages} and \Cref{eqn: cross messages}. Substituting \Cref{eqn: ttns node identity} into the expansion of \(\lVert P_{G} - \tilde{P} \rVert_{F}^{2}\) gives the two inner products
\begin{equation}\label{eqn: ttns inner products}
  \langle P_{G}, P_{G}\rangle = \Bigl\langle G_{q},\, \Bigl(\bigotimes_{v \in \mathcal{N}(q)} M_{v \to q}\Bigr) G_{q}\Bigr\rangle, \qquad
  \langle P_{G}, \tilde{P}\rangle = \Bigl\langle G_{q},\, \Bigl(\bigotimes_{v \in \mathcal{N}(q)} L_{v \to q}\Bigr) F_{q}\Bigr\rangle,
\end{equation}
and so the loss restricted to \(G_{q}\) is a convex quadratic plus the barrier, with gradient and Hessian
\begin{align}
  \nabla_{G_{q}} \ell &= 2 \Bigl(\bigotimes_{v \in \mathcal{N}(q)} M_{v \to q}\Bigr) G_{q} - 2 \Bigl(\bigotimes_{v \in \mathcal{N}(q)} L_{v \to q}\Bigr) F_{q} - \mu_{q}\, (1 \oslash G_{q}), \label{eqn: ttns grad}\\
  \nabla^{2}_{G_{q}} \ell &= 2 \Bigl(\bigotimes_{v \in \mathcal{N}(q)} M_{v \to q}\Bigr) + \mu_{q}\, \mathrm{diag}\left(1 \oslash (G_{q} \odot G_{q})\right), \label{eqn: ttns hess}
\end{align}
the analogues of \Cref{eqn: grad} and \Cref{eqn: hess}. \Cref{alg: fitting} runs unchanged: a forward and a backward sweep over the tree, one damped Newton step per visited component, and the adaptive rule \Cref{eqn: adaptive mu} with the sum over all entries of \(G_{q}\). The case split of \Cref{sec: NHT fitting extra} becomes the split between external and internal nodes. At an external node, the data Hessian acts on the bonds only and identically on each slice \(G_{q}(i, :)\), so the Newton solve decouples into \(n\) independent systems over the bond axes. At an internal node, we solve by the matrix-free CG of \Cref{sec: NHT fitting extra}, with one matrix product per incident edge.

\paragraph{Warm initialization}
The warm initialization needs only a notion of parent, and the prescribed root supplies it. One gauge \(Q_{(v,q)}\) of size \(\rho \times \rho\) sits on each edge, applied to the matching bond of both components that meet there. The two copies cancel in every contraction over the edge, as in \Cref{eqn: gauge cancel}, so the represented tensor does not depend on the gauges. As in \Cref{sec: NHT warm start}, we choose the gauges and the non-negative components so that each gauged component matches the signed \(F_{q}\). The decoupled sweep visits the nodes in leaf-to-root order of the prescribed root, so at \(q\) the gauges on the child edges are fixed and the parent-edge gauge \(Q\) is free. With \(Q_{q} := \bigotimes_{v \in \mathcal{N}(q)} Q_{(v,q)}\), and with the norm \(\lVert X \rVert_{M}^{2} := \langle X, (\bigotimes_{v \in \mathcal{N}(q)} M_{v \to q})\, X \rangle\) formed from the reference contractions \(\tilde{P}_{v \to q}\), the local fit at \(q\) is
\begin{equation}\label{eqn: ttns warm local}
  \min_{G_{q} > 0,\; Q^{\top} Q = I} \; \left\lVert Q_{q}\, G_{q} - F_{q} \right\rVert_{M}^{2} + \mu\, \phi_{q}(G_{q}),
\end{equation}
with the barrier coefficient \(\mu\) tracking the data term, as in \Cref{sec: NHT warm start}. In this norm the data term is exact: it equals the squared Frobenius error of \(\tilde{P}\) with node \(q\) alone replaced by its gauged fit. In diagram form, at an external node,
\begin{equation}\label{eqn: ttns warm diagram}
  \HTTTNSWarmLHS \;\approx\; \HTTTNSWarmRHS\,.
\end{equation}
The gauge step solves the Procrustes problem of \Cref{eqn: procrustes expand} matricized on the parent bond, with each fixed-gauge axis weighted by the Cholesky factor of its message. The component step is one interior-point Newton step on \Cref{eqn: ttns warm local}, with the Hessian \Cref{eqn: ttns hess} conjugated by the gauges. At the root, no parent edge exists, and so no gauge is solved.

\bibliographystyle{siamplain}
\bibliography{references}

\end{document}

%% file: ex_shared.tex
\usepackage{lipsum}
\usepackage{amsfonts}
\usepackage{graphicx}
\usepackage{epstopdf}
\usepackage{algorithmic}
\ifpdf
  \DeclareGraphicsExtensions{.eps,.pdf,.png,.jpg}
\else
  \DeclareGraphicsExtensions{.eps}
\fi

\headers{Density compression with non-negative hierarchical tucker}{X. Tang, H. Chen, L. Ying}

\title{Variational inference and density estimation with non-negative tensor of hierarchical tucker format}

\author{Xun Tang\thanks{Correspondence author, Department of Mathematics, Stanford University, Stanford, CA 94305, USA. 
  (\email{xuntang@stanford.edu}).}
  \and
  Haoxuan Chen\thanks{Institute for Computational and Mathematical Engineering (ICME), Stanford University, Stanford, CA 94305, USA. 
  (\email{haoxuanc@stanford.edu}).}
  \and
  Lexing Ying\thanks{Department of Mathematics and Institute for Computational and Mathematical Engineering (ICME), Stanford University, Stanford, CA 94305, USA. 
  (\email{lexing@stanford.edu})}
  \funding{X.T., H. C., and L.Y. are supported by AFOSR MURI award FA9550-24-1-0254.}
  }

\usepackage{amsopn}

\makeatletter
\newcommand*{\addFileDependency}[1]{%
  \typeout{(#1)}%
  \@addtofilelist{#1}%
  \IfFileExists{#1}{}{\typeout{No file #1.}}%
}
\makeatother

%% file: math_commands.tex
\usepackage{amsmath,amsfonts,bm}

\def\eqref#1{equation~\ref{#1}}

\def\1{\bm{1}}

\DeclareMathAlphabet{\mathsfit}{\encodingdefault}{\sfdefault}{m}{sl}
\SetMathAlphabet{\mathsfit}{bold}{\encodingdefault}{\sfdefault}{bx}{n}

\newcommand{\R}{\mathbb{R}}

%% file: ht_diagram.tex
\usetikzlibrary{positioning,calc}

\tikzset{
  tensor/.style={draw,rectangle,fill=blue!10,minimum width=1.8cm,
                 minimum height=1.1cm,thick,align=center},
  tensorv2/.style={draw,rounded corners=2pt,fill=blue!10,minimum width=1.8cm,
                 minimum height=1.1cm,thick,align=center},
  gauge/.style={draw,rounded corners=2pt,fill=blue!10,minimum width=0.9cm,
                 minimum height=0.9cm,thick,align=center},
  leg/.style={thick},
}
\newlength{\HTvleg}
\newlength{\HThleg}

\newcommand{\HTscale}{0.5}

\newcommand{\htinline}[2][1]{%
  \vcenter{\hbox{%
    \tikz[baseline={(current bounding box.center)},scale={#1*\HTscale},
          transform shape,every node/.style={font=\Large}]{#2}%
  }}%
}

\tikzset{
  htnode/.style={draw,rectangle,fill=blue!10,thick,
                 minimum size=0.95cm,inner sep=1pt,align=center},
  treeedge/.style={semithick},
}

\newcommand{\HTFigure}{%
  \def\D{-5}%
  \begin{tikzpicture}[x=1.5cm,y=1.2cm,every node/.style={font=\large}]
    \node[htnode] (r)  at (3.5,3) {$G^{(0)}_{1}$};
    \node[htnode] (a1) at (1.5,2) {$G^{(1)}_{1}$};
    \node[htnode] (a2) at (5.5,2) {$G^{(1)}_{2}$};
    \node[htnode] (b1) at (0.5,1) {$G^{(2)}_{1}$};
    \node[htnode] (b2) at (2.5,1) {$G^{(2)}_{2}$};
    \node[htnode] (b3) at (4.5,1) {$G^{(2)}_{3}$};
    \node[htnode] (b4) at (6.5,1) {$G^{(2)}_{4}$};
    \foreach \i in {1,...,8}{
      \node[htnode] (l\i) at ({\i-1},0) {$G^{(3)}_{\i}$};
    }
    \draw[leg] (r)--(a1) (r)--(a2);
    \draw[leg] (a1)--(b1) (a1)--(b2) (a2)--(b3) (a2)--(b4);
    \draw[leg] (b1)--(l1) (b1)--(l2) (b2)--(l3) (b2)--(l4)
               (b3)--(l5) (b3)--(l6) (b4)--(l7) (b4)--(l8);
    \foreach \i in {1,...,8}{
      \draw[leg] (l\i.south)--++(0,-0.55);
      \node[below] at ($(l\i.south)+(0,-0.55)$) {$i_{\i}$};
    }
    \node (tr)  at (3.5,{3+\D}) {$\{1,\dots,8\}$};
    \node (ta1) at (1.5,{2+\D}) {$\{1,2,3,4\}$};
    \node (ta2) at (5.5,{2+\D}) {$\{5,6,7,8\}$};
    \node (tb1) at (0.5,{1+\D}) {$\{1,2\}$};
    \node (tb2) at (2.5,{1+\D}) {$\{3,4\}$};
    \node (tb3) at (4.5,{1+\D}) {$\{5,6\}$};
    \node (tb4) at (6.5,{1+\D}) {$\{7,8\}$};
    \foreach \i in {1,...,8}{
      \node (tl\i) at ({\i-1},{0+\D}) {$\{\i\}$};
    }
    \draw[treeedge] (tr)--(ta1) (tr)--(ta2);
    \draw[treeedge] (ta1)--(tb1) (ta1)--(tb2) (ta2)--(tb3) (ta2)--(tb4);
    \draw[treeedge] (tb1)--(tl1) (tb1)--(tl2) (tb2)--(tl3) (tb2)--(tl4)
                    (tb3)--(tl5) (tb3)--(tl6) (tb4)--(tl7) (tb4)--(tl8);
  \end{tikzpicture}%
}

\newcommand{\HTFull}[2][1]{%
  \htinline[#1]{%
    \node[tensor,minimum width=3.6cm] (C) {$#2$};
    \foreach \x/\lab in {0.6/{i_a}, 1.8/{i_b}, 3.0/{i_f}}{
      \draw[leg] ($(C.south west)+(\x,0)$)--++(0,-\HTvleg);
      \node[below] at ($(C.south west)+(\x,-\HTvleg)$) {$\lab$};
    }
  }%
}

\newcommand{\HTFullBi}[2][1]{%
  \htinline[#1]{%
    \node[tensor,minimum width=3.6cm] (C) {$#2$};
    \foreach \x/\lab in {1.2/{i_{q}}, 2.4/{i_{f}}}{
      \draw[leg] ($(C.south west)+(\x,0)$)--++(0,-\HTvleg);
      \node[below] at ($(C.south west)+(\x,-\HTvleg)$) {$\lab$};
    }
  }%
}

\newcommand{\HTUnfold}[1][1]{%
  \htinline[#1]{%
    \node[tensor,minimum width=2.0cm] (Cbar) {$C_{q}$};
    \node[tensor,minimum width=2.0cm,right=1.3cm of Cbar] (Cf) {$C_{f}$};
    \draw[leg] (Cbar.east)--(Cf.west);
    \draw[leg] (Cbar.south)--++(0,-\HTvleg);
    \node[below] at ($(Cbar.south)+(0,-\HTvleg)$) {$i_{q}$};
    \draw[leg] (Cf.south)--++(0,-\HTvleg);
    \node[below] at ($(Cf.south)+(0,-\HTvleg)$) {$i_{f}$};
  }%
}

\newcommand{\HTNodeMap}[5][1]{%
  \htinline[#1]{%
    \node[tensor,minimum width=1.8cm] (Ca) at (0,0) {$#3$};
    \node[tensor,minimum width=1.8cm] (Cb) at (2.9,0) {$#4$};
    \node[tensor,minimum width=1.8cm] (Cf) at (5.8,0) {$#5$};
    \node[tensor,minimum width=1.4cm] (Fq) at (2.9,1.9) {$#2$};
    \draw[leg] (Fq.south) to (Ca.north);
    \draw[leg] (Fq.south) to (Cb.north);
    \draw[leg] (Fq.south) to (Cf.north);
    \foreach \B/\lab in {Ca/{i_a}, Cb/{i_b}, Cf/{i_f}}{
      \draw[leg] (\B.south)--++(0,-\HTvleg);
      \node[below] at ($(\B.south)+(0,-\HTvleg)$) {$\lab$};
    }
  }%
}

\newcommand{\HTSolveLHS}[1][1]{%
  \htinline[#1]{%
    \node[tensor,minimum width=1.5cm] (Aa) at (0,0)   {$A_{a}$};
    \node[tensor,minimum width=1.5cm] (Ab) at (2.6,0) {$A_{b}$};
    \node[tensor,minimum width=1.5cm] (Af) at (5.2,0) {$A_{f}$};
    \node[tensor,minimum width=1.4cm] (Fq) at (2.6,1.9) {$F_q$};
    \draw[leg] (Fq.south) to (Aa.north);
    \draw[leg] (Fq.south) to (Ab.north);
    \draw[leg] (Fq.south) to (Af.north);
    \foreach \B/\lab in {Aa/{\alpha}, Ab/{\beta}, Af/{\gamma}}{
      \draw[leg] (\B.south)--++(0,-\HTvleg);
      \node[below] at ($(\B.south)+(0,-\HTvleg)$) {$\lab$};
    }
  }%
}
\newcommand{\HTSolveRHS}[1][1]{%
  \htinline[#1]{%
    \node[tensor,minimum width=2.2cm] (Bq) {$B_q$};
    \foreach \x/\lab in {0.5/{\alpha},1.1/{\beta},1.7/{\gamma}}{
      \draw[leg] ($(Bq.south west)+(\x,0)$)--++(0,-\HTvleg);
      \node[below] at ($(Bq.south west)+(\x,-\HTvleg)$) {$\lab$};
    }
  }%
}

\tikzset{
  gleg/.style={line width=2pt},
}

\newcommand{\HTTTNSFigure}{%
  \begin{tikzpicture}[x=1.5cm,y=1.2cm,every node/.style={font=\large}]
    \node[htnode] (r)  at (2.5,2) {$G_{r}$};
    \node[htnode] (x1) at (0.5,1) {$G_{1}$};
    \node[htnode] (w)  at (2.0,1) {$G_{w}$};
    \node[htnode] (x2) at (4.2,1) {$G_{2}$};
    \node[htnode] (x3) at (1.3,0) {$G_{3}$};
    \node[htnode] (x4) at (2.7,0) {$G_{4}$};
    \node[htnode] (x5) at (4.2,0) {$G_{5}$};
    \draw[leg] (r)--(x1) (r)--(w) (r)--(x2) (w)--(x3) (w)--(x4) (x2)--(x5);
    \foreach \B/\lab in {x1/{i_{1}}, x3/{i_{3}}, x4/{i_{4}}, x5/{i_{5}}}{
      \draw[leg] (\B.south)--++(0,-0.55);
      \node[below] at ($(\B.south)+(0,-0.55)$) {$\lab$};
    }
    \draw[leg] ($(x2.south)+(0.3,0)$)--++(0.3,-0.55);
    \node[below] at ($(x2.south)+(0.6,-0.55)$) {$i_{2}$};
    \node[above=3pt of r] {root};
  \end{tikzpicture}%
}

\newcommand{\HTTTNSSolveLHS}[1][1]{%
  \htinline[#1]{%
    \node[tensor,minimum width=1.5cm] (Fq) at (0,0) {$F_{q}$};
    \node[tensor,minimum width=1.5cm] (A) at (2.4,0) {$A_{q}$};
    \draw[leg] (Fq.east)--(A.west);
    \draw[gleg] (A.east)--++(\HThleg,0);
    \node[right] at ($(A.east)+(\HThleg,0)$) {$\alpha_{\mathcal{E}(q)}$};
    \draw[leg] (Fq.south)--++(0,-\HTvleg);
    \node[below] at ($(Fq.south)+(0,-\HTvleg)$) {$i_{q}$};
  }%
}
\newcommand{\HTTTNSSolveRHS}[1][1]{%
  \htinline[#1]{%
    \node[tensor,minimum width=1.5cm] (Bq) {$B_{q}$};
    \draw[gleg] (Bq.east)--++(\HThleg,0);
    \node[right] at ($(Bq.east)+(\HThleg,0)$) {$\alpha_{\mathcal{E}(q)}$};
    \draw[leg] (Bq.south)--++(0,-\HTvleg);
    \node[below] at ($(Bq.south)+(0,-\HTvleg)$) {$i_{q}$};
  }%
}

\newcommand{\HTTTNSWarmLHS}[1][1]{%
  \htinline[#1]{%
    \node[tensor,minimum width=1.5cm] (Cq) at (0,0) {$C_{q}$};
    \node[gauge] (Qq) at (2.2,0) {$Q_{q}$};
    \draw[leg] (Cq.east)--(Qq.west);
    \draw[gleg] (Qq.east)--++(\HThleg,0);
    \node[right] at ($(Qq.east)+(\HThleg,0)$) {$\alpha_{\mathcal{E}(q)}$};
    \draw[leg] (Cq.south)--++(0,-\HTvleg);
    \node[below] at ($(Cq.south)+(0,-\HTvleg)$) {$i_{q}$};
  }%
}
\newcommand{\HTTTNSWarmRHS}[1][1]{%
  \htinline[#1]{%
    \node[tensor,minimum width=1.5cm] (Fq) {$F_{q}$};
    \draw[gleg] (Fq.east)--++(\HThleg,0);
    \node[right] at ($(Fq.east)+(\HThleg,0)$) {$\alpha_{\mathcal{E}(q)}$};
    \draw[leg] (Fq.south)--++(0,-\HTvleg);
    \node[below] at ($(Fq.south)+(0,-\HTvleg)$) {$i_{q}$};
  }%
}

\newcommand{\HTGaugeEdge}[1][1]{%
  \htinline[#1]{%
    \node[tensor,minimum width=1.5cm] (Cp) at (0,3.4) {$G_{f}$};
    \node[gauge] (gA) at (0,2.1) {$Q_{e}$};
    \node[gauge] (gB) at (0,0.9) {$Q_{e}$};
    \node[tensor,minimum width=1.5cm] (Cq) at (0,-0.4) {$G_{q}$};
    \draw[leg] (Cp.south)--(gA.north);
    \draw[leg] (gA.south)--(gB.north);
    \draw[leg] (gB.south)--(Cq.north);
    \draw[leg] (Cp.north)--++(0,\HTvleg);
    \draw[leg] ($(Cp.south)+(-0.55,0)$)--++(-0.25,-\HTvleg);
    \draw[leg] ($(Cq.south)+(-0.35,0)$)--++(0,-\HTvleg);
    \draw[leg] ($(Cq.south)+(0.35,0)$)--++(0,-\HTvleg);
  }%
}
\newcommand{\HTGaugeEdgeRHS}[1][1]{%
  \htinline[#1]{%
    \node[tensor,minimum width=1.5cm] (Cp) at (0,1.4) {$G_{f}$};
    \node[tensor,minimum width=1.5cm] (Cq) at (0,0) {$G_{q}$};
    \draw[leg] (Cp.south)--(Cq.north);
    \draw[leg] (Cp.north)--++(0,\HTvleg);
    \draw[leg] ($(Cp.south)+(-0.55,0)$)--++(-0.25,-\HTvleg);
    \draw[leg] ($(Cq.south)+(-0.35,0)$)--++(0,-\HTvleg);
    \draw[leg] ($(Cq.south)+(0.35,0)$)--++(0,-\HTvleg);
  }%
}

\newcommand{\HTWarmLHS}[1][1]{%
  \htinline[#1]{%
    \node[gauge] (Qa) at (0,0)   {$Q_{a}$};
    \node[gauge] (Qb) at (2.6,0) {$Q_{b}$};
    \node[gauge] (Qp) at (5.2,0) {$Q$};
    \node[tensor,minimum width=1.4cm] (Cq) at (2.6,1.9) {$G_{q}$};
    \draw[leg] (Cq.south) to (Qa.north);
    \draw[leg] (Cq.south) to (Qb.north);
    \draw[leg] (Cq.south) to (Qp.north);
    \foreach \B/\lab in {Qa/{\alpha}, Qb/{\beta}, Qp/{\gamma}}{
      \draw[leg] (\B.south)--++(0,-\HTvleg);
      \node[below] at ($(\B.south)+(0,-\HTvleg)$) {$\lab$};
    }
  }%
}
\newcommand{\HTWarmRHS}[1][1]{%
  \htinline[#1]{%
    \node[tensor,minimum width=2.2cm] (Fq) {$F_q$};
    \foreach \x/\lab in {0.5/{\alpha},1.1/{\beta},1.7/{\gamma}}{
      \draw[leg] ($(Fq.south west)+(\x,0)$)--++(0,-\HTvleg);
      \node[below] at ($(Fq.south west)+(\x,-\HTvleg)$) {$\lab$};
    }
  }%
}

\newcommand{\HTProcK}[1][1]{%
  \htinline[#1]{%
    \node[tensor,minimum width=1.5cm] (F) {$\hat{F}_{q}$};
    \node[tensor,minimum width=1.5cm,below=0.9cm of F] (C) {$\hat{C}_{q}$};
    \node[tensor,minimum width=1.0cm,right=0.6cm of F] (M) {$M_{f}$};
    \draw[leg] ($(F.south west)+(0.45,0)$)--($(C.north west)+(0.45,0)$);
    \draw[leg] ($(F.south east)+(-0.45,0)$)--($(C.north east)+(-0.45,0)$);
    \draw[leg] (F.east)--(M.west);
    \draw[leg] (M.east)--++(\HThleg,0);
    \node[right] at ($(M.east)+(\HThleg,0)$) {$\gamma$};
    \draw[leg] (C.east)--++(\HThleg,0);
    \node[right] at ($(C.east)+(\HThleg,0)$) {$\gamma'$};
  }%
}

\newcommand{\HTGram}[3][1]{%
  \htinline[#1]{%
    \node[tensor,minimum width=1.8cm] (top) {$#2$};
    \node[tensor,minimum width=1.8cm,below=0.9cm of top] (bot) {$#2$};
    \draw[leg] (top.south)--(bot.north);
    \draw[leg] (top.east)--++(\HThleg,0);
    \node[right] at ($(top.east)+(\HThleg,0)$) {$\alpha$};
    \draw[leg] (bot.east)--++(\HThleg,0);
    \node[right] at ($(bot.east)+(\HThleg,0)$) {$\alpha'$};
  }%
}

\newcommand{\HTCross}[1][1]{%
  \htinline[#1]{%
    \node[tensor,minimum width=1.8cm] (top) {$P_{G,a}$};
    \node[tensor,minimum width=1.8cm,below=0.9cm of top] (bot) {$\tilde{P}_{a}$};
    \draw[leg] (top.south)--(bot.north);
    \draw[leg] (top.east)--++(\HThleg,0);
    \node[right] at ($(top.east)+(\HThleg,0)$) {$\alpha$};
    \draw[leg] (bot.east)--++(\HThleg,0);
    \node[right] at ($(bot.east)+(\HThleg,0)$) {$\alpha'$};
  }%
}